\documentclass{article}
\usepackage[utf8]{inputenc}
\usepackage{amsmath}
\usepackage{algorithm}
\usepackage{algpseudocode}
\usepackage{amsthm}
\usepackage{amsfonts}
\usepackage{amssymb} 
\usepackage{amssymb}
\usepackage{hyperref}
\usepackage{xcolor}
\usepackage{braket}
\hypersetup{
    colorlinks,
    linkcolor={red!80!black},
    citecolor={blue!80!black},
    urlcolor={blue!80!black}
}
\usepackage{float}
\usepackage[margin=1in]{geometry}
\usepackage{graphicx}
\usepackage{subcaption}
\usepackage{caption}
\usepackage{url}
\usepackage{bm}
\usepackage{enumitem}
\usepackage{mathtools}
\usepackage{authblk}

\newtheorem{remark}{Remark}

\newtheorem{theorem}{Theorem}
\newtheorem{proposition}{Proposition}

\newcommand{\Real}{\mathbb R}
\newcommand{\Complex}{\mathbb C}

\newcommand{\qbar}{\overline{q}}

\newcommand{\us}{u^{*}}

\newcommand{\bfqone}{\mathbf{q}_1}
\newcommand{\bfq}{\mathbf{q}}
\newcommand{\bfzone}{\mathbf{z}_1}
\newcommand{\bfwone}{\mathbf{w}_1}

\newcommand{\bfzero}{\mathbf{0}}

\newcommand{\bfk}{\mathbf{k}}
\newcommand{\bfn}{\mathbf{n}}

\newcommand{\bfqtone}{\tilde{\mathbf{q}}_1}
\newcommand{\qtzero}{\tilde{q}_{0}}

\newcommand{\Hmod}{\bar{H}}

\newcommand{\bfx}{\mathbf{x}}

\newcommand{\qt}{\widetilde{q}}
\newcommand{\bq}{\mathbf{q}}
\newcommand{\bR}{\mathbf{R}}
\newcommand{\bqt}{\widetilde{\mathbf{q}}}

\newcommand{\eps}{\varepsilon}
\newcommand{\bfdelta}{\boldsymbol{\delta}}

\newcommand{\ktk}{|\bfk|^2}
\newcommand{\slowmani}{\mathcal{K}_{\tau}}
\newcommand{\slowmap}{\mathcal{M}_{\tau}[q_0]}
\newcommand{\smapiter}[1]{\mathcal{M}_{\tau}^{#1}[q_0]}
\newcommand{\induced}{\mathcal{N}_{\tau}[q_0]}
\newcommand{\inditer}[1]{\mathcal{N}_{\tau}^{#1}[q_0]}

\newcommand{\bigo}{\mathcal{O}}
\newcommand{\dx}{{\Delta x}}
\newcommand{\dy}{{\Delta y}}
\newcommand{\dt}{{\Delta t}}
\newcommand{\Courant}{{\mathcal{C}}}
\newcommand{\refcrit}{{\Delta \rho_\text{ref}}}

\DeclareMathOperator{\diag}{diag}

\title{A hyperbolic nonlinear Schrödinger system
in multiple spatial dimensions: structure, convergence, and numerical discretization}
\author[1]{Carlos Mu\~noz-Moncayo \footnote{Corresponding author. \url{carlos.munozmoncayo@kaust.edu.sa}}}
\author[2]{Abhijit Biswas}
\author[1]{David I. Ketcheson} 
\author[1]{Manvendra Rajvanshi}

\affil[1]{Computer Electrical and Mathematical Science and Engineering Division
(CEMSE), King Abdullah University of Science and Technology (KAUST),
Thuwal, 23955-6900, Saudi Arabia}
\affil[2]{Indian Institute of Technology Kanpur, Kalyanpur, Kanpur, 208016, Uttar Pradesh, India}

\begin{document}

\maketitle
\begin{abstract}
We propose and study a first-order hyperbolic approximation of the nonlinear Schrödinger (NLS) equation in $N$ spatial dimensions.
This is an extension of the recently-proposed 1D hyperbolic NLS model.
We study the structure and properties of the hyperbolic approximation and their
relation to those of the original system.  We show that the N-dimensional
hyperbolic system exactly conserves analogs of the mass, energy, and momentum,
and that it possesses stationary solutions that are closely related to those of NLS.
The hyperbolic system is also shown to admit
plane-wave solutions that share some stability properties with those of NLS,
and that relaxation can shrink the focusing low-wavenumber modulational instability band
as well as introduce a high-wavenumber instability band not present in NLS.
We also show that solutions of the hyperbolic system converge to solutions of NLS
in the relaxation limit.
We present an adaptive-mesh finite volume discretization of the hyperbolic system
and demonstrate its effectiveness in handling simulations, including those
with geometric obstacles or non-periodic boundary conditions.  The effectiveness
of the numerical scheme is demonstrated through examples including a Bose-Einstein
condensate, line soliton instability, and superflow past a cylinder.
\end{abstract}

\section{Introduction}
In this work we propose and study the following system in $N$ spatial dimensions:
\begin{subequations}
    \label{eq:ND-NLSH}
\begin{align}
     \label{eq:ND-NLSH q0}
    i \partial_t q_0 + \mathrm{div}(\mathbf{q}_1) + (V_0(\bfx)+\kappa |q_0|^{p-1}) q_0 & = 0,\\
    \label{eq:ND-NLSH q1}
   i \tau \partial_t \mathbf{q}_1 - \nabla q_0 + \mathbf{q}_1 & = 0,
\end{align}
\end{subequations}
where $\mathbf{q}_1 = (q_1, q_2, \ldots, q_N)^T$ and each $q_j(\bfx,t)$ is
complex-valued.  This is a first-order hyperbolic system, but due to the
presence of the algebraic term in the second equation it exhibits dispersive
behavior.  In fact, as we will see it possesses solutions that are typical of
higher-order dispersive nonlinear wave equations, including traveling solitary
waves, dispersive shock waves, and transverse and modulational instabilities.
We refer to it as a hyperbolic nonlinear Schrödinger system (NLSH for short)
because, in the limit when the relaxation parameter $\tau$ goes to zero, it
formally becomes equivalent to the nonlinear Schrödinger equation (NLS):
\begin{align}\label{eq:ND-NLS}
    i \partial_t u + \nabla^2 u + (V_0(\bfx)+\kappa |u|^{p-1}) u = 0.
\end{align}
NLS is a canonical model for a variety of physical problems, including the dynamics
of Bose-Einstein condensates, surface water waves, and nonlinear optics \cite{sulem2007nonlinear}.
Here $u: \Real^N \times \Real \to \mathbb{C}$ is a complex-valued function, $\kappa \in \Real\setminus\{0\}$ is a constant,
and $V_0: \Real^N \to \Real$ is a potential function.
As mentioned already, solutions to \eqref{eq:ND-NLS} display a rich variety of structures, including solitons,
transverse and modulational instabilities, and
dispersive shockwaves.
The accurate numerical simulation of these phenomena has been the subject 
of extensive research over the past several decades, and remains an active area of study.

Like NLS, traditional models of dispersive nonlinear waves involve higher-order PDEs.
Dispersive first-order hyperbolic systems can be a useful and interesting alternative
for both analytical and computational reasons. Hyperbolic systems come with a
well-developed theory based on characteristics, finite propagation speeds,
symmetrization, and energy estimates, and these structural features may reveal
properties that are less transparent in high-order equations. From the
numerical point of view, first-order hyperbolic systems are amenable to a broad
range of local discretizations, including finite volume and discontinuous
Galerkin methods, and are particularly well suited to adaptive mesh refinement,
complex geometries, and nonperiodic boundary conditions. When such a
reformulation contains a relaxation parameter and converges to the original
equation in an appropriate singular limit, it also provides a controlled way of
trading modeling error against the structural and computational advantages of
the hyperbolic formulation.

Two hyperbolic approximations to the one-dimensional NLS equation have recently been developed.
In the first one, the hydrodynamic formulation of the equation
is augmented into a hyperbolic system through an extended Lagrangian approach \cite{dhaouadi2019extended}.
The system is further generalized to the higher dimensional setting and discretized with
high-order ADER schemes in \cite{busto2021high}.
This approach is
naturally limited to the defocusing regime and degenerates at
vacuum points, where the Madelung transform is singular.
In the second approach \cite{ketcheson2025approximation,biswas2025hyperbolic},
the relaxation is done directly at the level of the wavefunction
by means of an auxiliary variable that approximates the gradient of the solution.
The resulting system is is known to possess conserved mass, momentum,
and energy functionals analogous to those of NLS,
and is equally valid in the focusing and defocusing regimes.

The NLSH system \eqref{eq:ND-NLSH} is a generalization of the latter approach.
It is interesting in that it not only approximates NLS in the relaxation limit,
but it possesses remarkable properties, many related to those of NLS, for arbitrary
values of $\tau$ far from the limit. 
Our purpose is to develop the mathematical theory of this hyperbolic
Schrödinger system, establish rigorously its relationship with NLS, and
investigate the numerical opportunities created by its first-order hyperbolic
structure.  

We first investigate the properties of NLSH solutions, focusing on the similarities
and differences relative to those of NLS.  In Section \ref{sec:basics}, we derive a Hamiltonian structure and
conserved mass, momentum, and energy.  In Section \ref{sec:special}, we show that
NLSH admits plane wave and standing wave solutions that are analogous to those of NLS,
and that converge to the latter uniformly and linearly in $\tau$.  We also
perform a modulational stability analysis, which reveals that hyperbolization can shrink the focusing low-wavenumber
instability band of NLS and introduce a high-wavenumber unstable band with no
NLS counterpart. The latter is suppresed for every wavenumber provided
that $\tau$ lies below a threshold depending on the nonlinearity strength.
In Section \ref{sec:convergence}, we prove (following the approach of
Giesselmann \& Ranocha \cite{giesselmann2025convergence})
that sufficiently smooth solutions of NLSH converge to those
of NLS, thus providing the first rigorous justification for hyperbolic approximation of NLS.
Initializations of higher formal order in $\tau$ of the auxiliary variables are
characterized  through a formal slow-manifold analysis in Section \ref{sec: well-prepared initial data}.

Turning to the numerical discretization of NLS,
perhaps the most popular method for solving \eqref{eq:ND-NLS} is the Fourier spectral
method with operator splitting \cite{bao2003numerical}, in which one alternates between evolving the dynamics
given by the spatial derivative term and those given by the algebraic terms.
This method can be extremely efficient and accurate, but it
is mostly suited to uniform grids in very simple geometries
with periodic boundary conditions. 
For more complicated settings, common alternatives follow the method of lines,
discretizing first in space with finite differences \cite{Delfour1981FinitedifferenceSO},
finite elements \cite{katsaounis2018posteriori}, or
local discontinuous Galerkin methods \cite{xu2005local};
see \cite{antoine2013computational} for a comprehensive review.
Because
the Laplacian term induces stiffness, such methods require the use of implicit
time discretization and the solution of a linear algebraic system at each step.

For moderate values of $\tau$, NLSH is not stiff, making it amenable to explicit discretizations.
A potential advantage of discretizing NLSH instead of NLS is therefore
the possibility to avoid the need for implicit time stepping.  
The advantage of explicit integration is more pronounced
in higher spatial dimensions since the cost of the linear algebraic solve increases
dramatically.  Additionally, first-order hyperbolic systems are well-suited to
the use of adaptive mesh refinement (AMR) in order to resolve small-scale features
that may appear in localized regions.  The first-order approach allows one to take
advantage the well-developed ecosystem of numerical methods for hyperbolic PDEs
and to impose boundary conditions in a well-studied framework. 
In Section \ref{sec:FVM}, we present an exemplar of this family of methods: an explicit
finite volume discretization of NLSH that is compatible with mapped grids, embedded boundaries, and
patch-based AMR.
The solver is tested against a variety of complex scenarios, including gray soliton reflections, transverse instabilities
of bright solitons, radially symmetric dispersive shock waves, and vortex shedding past a cylinder
in Sections \ref{sec:soliton_reflections}-\ref{sec:flow past cylinder}.
Furthermore, the theoretical results of Sections \ref{sec:basics}-\ref{sec:convergence} are verified numerically
in Sections \ref{sec:MI numerical}-\ref{sec:numerical_ground_state} with an asymptotic-preserving implicit-explicit (ImEx) scheme
on tests with periodic uniform grids.

NLSH provides an alternative route to NLS computation in which a quantitatively
controlled modeling error is exchanged for a first-order hyperbolic structure,
based on the choice of $\tau$.
A central contribution of this work is to
quantify the former and, in particular, provide practical estimates
such that the choice of relaxation parameter is more than an heuristic.
We do not
expect the first-order approach to be advantageous in all situations;
nevertheless, taken together, these potential benefits seem compelling
and worthy of exploration.

The numerical results presented in this work can be reproduced 
using the code in the repository \cite{munoz2026repository}.

\section{The First-order Hyperbolic System and its Hamiltonian Structure} \label{sec:basics}

We can write \eqref{eq:ND-NLSH} in the canonical form
\begin{align} \label{hyperbolic-form}
    A_0 \partial_t \begin{bmatrix}
        q_0\\
        \mathbf{q}_1
    \end{bmatrix} +
    \sum_{j=1}^{N} A_j \partial_{j} \begin{bmatrix}
        q_0\\
        \mathbf{q}_1
    \end{bmatrix}
        = \begin{bmatrix}
            i\kappa |q_0|^{p-1} q_0 + i V_0(\bfx)q_0 \\
            i \mathbf{q}_1
        \end{bmatrix},
\end{align}
where $A_0 = \diag(1,\tau, \dots, \tau)$ and
\begin{align*}
    A_j = \begin{bmatrix}
        0 & -i \mathbf{e}_j^T\\
        i \mathbf{e}_j & \mathbf{0}_{N\times N}
    \end{bmatrix} \; ,
\end{align*}
and $\mathbf{e}_j$ is the $j$-th unit vector in $\mathbb{R}^N$.
It is indeed hyperbolic (see e.g. \cite{leveque2002finite}), since the matrix
\begin{align*}
    A = A_0^{-1}\sum_{j=1}^{N} a_j A_j 
      = \begin{bmatrix}
        0 & -i \mathbf{a}^T \\
        i \tau^{-1}\mathbf{a} & \mathbf{0}_{N\times N}
    \end{bmatrix} ,
\end{align*}
has real eigenvalues
$
   \lambda \in \{\pm \tau^{-1/2}\|\mathbf{a}\|_2, 0\}$.

Before deriving the Hamiltonian structure of NLSH, we recall the corresponding
structure of NLS. 
For simplicity, throughout this and the next section we will assume that $V_0\equiv 0$.
\subsection{Hamiltonian Structure of NLS}
Let $u^*$ denote the complex conjugate of $u$. Then the $N$-dimensional NLS equation \eqref{eq:ND-NLS} can be written as 
\begin{align*}
    iu_t & = \frac{\delta H}{\delta u^*} \;,
\end{align*}
where
\begin{align}\label{Eq:Hamilton_NLS}
H = \int_{\Real^N} \left(|\nabla u|^2 - \frac{2\kappa}{p+1}|u|^{p+1}\right) dx,
\end{align}
and $\frac{\delta H}{\delta \us}$ is the variational derivative of the Hamiltonian with respect to $\us$. For any solution of the NLS equation, the mass,
momentum, and energy are conserved. They are given by
\begin{subequations} \label{NLS-conserved}
\begin{align}
    M[u(t)] &= \int_{\Real^N} |u(x,t)|^2 \, dx \;,\\
    P_j[u(t)] &= \int_{\Real^N} \mathrm{Im}(u^* \partial_j u) \, dx \;, \qquad  \forall j=1,\ldots,N,\\
    E[u(t)] &= \int_{\Real^N} \left(|\nabla u(x,t)|^2 - \frac{2\kappa}{p+1}|u(x,t)|^{p+1}\right) dx \;. \label{eq:energy NLS}
\end{align}
\end{subequations}
\subsection{Hamiltonian Structure of NLSH}
The $N$-dimensional NLSH system of equations possesses an equilibrium structure similar to that of the NLS equation. We consider the modified Hamiltonian $\Hmod$:
\begin{align}
    \Hmod = \int_{\Real^N} \left( \mathbf{q}_{1}^{*} \cdot \nabla q_{0}  + \mathbf{q}_{1} \cdot \nabla q_{0}^{*} - |\mathbf{q}_1|^2 - \frac{2\kappa}{p+1}|q_0|^{p+1}\right) dx \;,
    \label{Eq:Hamilton_NLSH_modified}
\end{align} 
where $q_0^{*}$ and $\mathbf{q}_1^{*}$ denote the complex conjugates of the fields $q_0$ and $\mathbf{q}_1$, respectively. The $N$-dimensional NLSH system \eqref{eq:ND-NLSH} can be written in Hamiltonian PDE form as
\begin{align}
    \label{eq:NLSH as a Hamiltonian}
    i Q_t = \bar{\mathcal{J}} \frac{\delta \bar{H}}{\delta Q^*},
\end{align}
with
\begin{align}
    Q =  \begin{bmatrix}
         q_0 \\
         \mathbf{q}_1 
     \end{bmatrix}, \quad
     \bar{\mathcal{J}} = 
     \begin{bmatrix}
        1 & \mathbf{0}^{T} \\
        \mathbf{0} & \tau^{-1}\mathbf{I}_{N \times N} 
     \end{bmatrix}, \quad
     \frac{\delta \bar{H}}{\delta Q^*} = 
     \begin{bmatrix}
        \frac{\delta \bar{H}}{\delta q_0^{*}} \\
        \frac{\delta \bar{H}}{\delta \mathbf{q}_1^{*}}   
     \end{bmatrix} \;,
\end{align}
where $\mathbf{0} \in \Real^{N}$. The variational derivatives $\frac{\delta \bar{H}}{\delta q_0^{*}}$ and $\frac{\delta \bar{H}}{\delta \mathbf{q}_1^{*}}$ of the modified Hamiltonian $\Hmod$ with respect to $q_0^{*}$ and $\mathbf{q}_1^{*}$, respectively, are given by
\begin{align*}
     \begin{bmatrix}
        \frac{\delta \bar{H}}{\delta q_0^{*}} \\
        \frac{\delta \bar{H}}{\delta \mathbf{q}_1^{*}} 
     \end{bmatrix} & = 
     \begin{bmatrix}
         - \mathrm{div}(\mathbf{q}_1) - \kappa |q_0|^{p-1} q_0 \\
         \nabla q_0 - \mathbf{q}_1  
     \end{bmatrix}
     \;.
\end{align*}
It is straightforward to see that the modified NLSH Hamiltonian \eqref{Eq:Hamilton_NLSH_modified} is naturally related
to that of NLS \eqref{Eq:Hamilton_NLS} through the following relation
\begin{align*}
    \Hmod[Q(t)] &= \int_{\Real^N} \left( 2\mathrm{Re}(\bfqone^*\nabla q_0) - |\mathbf{q}_1|^2 - \frac{2\kappa}{p+1}|q_0|^{p+1}\right) dx \\
    &= \int_{\Real^N} \left(|\nabla q_0|^2 - |\nabla q_0-\mathbf{q}_1|^2 - \frac{2\kappa}{p+1}|q_0|^{p+1}\right) dx \\
    &= H[q_0(t)] -||\nabla q_0-\mathbf{q}_1||_{L^2(\Real^N)}^2.
\end{align*}
Similarly to the one-dimensional case, it is possible to show that the NLSH system shares several properties
with the NLS equation, including the existence of modified conserved quantities that converge to those of NLS in the relaxation limit.
This is summarized in the following proposition.

\begin{proposition}
    \label{prop: conserved quantities}
    For $V_0\equiv 0$ and $Q$ smooth and decaying fast enough away from the origin, 
    the NLSH system possesses three conserved quantities analogous to the mass, momentum, and energy of the NLS equation; 
    these are given respectively by
    \begin{subequations}
    \begin{align}
    \bar{M}[Q(t)] &= \int_{\Real^N} (|q_0|^2+\tau|\mathbf{q}_1|^2) \, dx \;,\label{eq:modified mass NLSH}\\
    \bar{P}_j[Q(t)] &= \int_{\Real^N} \mathrm{Im}(q_0^{*} \partial_j q_0) +
    \tau \sum_{k=1}^{N} \mathrm{Im}(q_k^{*} \partial_j q_k) \, dx\;, \qquad  \forall j=1,\ldots,N\;, \label{eq:modified momentum NLSH}\\
    \bar{H}[Q(t)] &= \int_{\Real^N} \left( \mathbf{q}_{1}^{*} \cdot \nabla q_{0}  +
    \mathbf{q}_{1} \cdot \nabla q_{0}^{*} - |\mathbf{q}_1|^2 - \frac{2\kappa}{p+1}|q_0|^{p+1}\right) dx \;.\label{eq:modified energy NLSH}
    \end{align}
    \end{subequations}
\end{proposition}
\begin{proof}
    \begin{enumerate}
        \item The conservation of the modified mass $\bar{M}$ follows after noting that
        \begin{align}
        \int_{\Real^N} \frac{\partial}{\partial t}\left[\left(|q_0|^2 + \tau |\mathbf{q}_1|^2 \right) \right] dx =
        i \int_{\Real^N} \mathrm{div} \left(q_0^{*} \mathbf{q}_1 - q_0 \mathbf{q}_1^{*} \right) dx =
        i \int_{\Real^N} \left[ \sum_{j=1}^{N} \frac{\partial}{\partial x_j}\left(q_0^{*}q_j-q_0q_j^{*}\right)  \right] dx\;.
        \end{align}
        \item For the conservation of the modified momentum $\bar{P}_j$, we first conjugate the NLSH system \eqref{eq:ND-NLSH} to get
            \begin{subequations}
                \label{eq:conjugated system}
            \begin{align}
                &\frac{\partial}{\partial t} q_0^{*} = - i \mathrm{div}\left(\mathbf{q}_1^{*}\right) - i \kappa |q_0|^{p-1} q_0^{*},\\
                &\tau \frac{\partial}{\partial t} q_k^{*} = i \partial_k q_0^{*} - i q_k^{*}, \quad k=1,\ldots,N,
            \end{align}
            \end{subequations}
            and observe that for every $k=0,\ldots,N$, we have
            \begin{align}
                \label{eq:claim imaginary 1}
                \frac{d}{dt} \mathrm{Im} \int_{\Real^N} q_k^{*} \frac{\partial}{\partial x_j} q_k \, d\bfx
                &= 2\mathrm{Im} \int_{\Real^N} \frac{\partial}{\partial t} q_k^{*} \frac{\partial}{\partial x_j} q_k \, d\bfx.
            \end{align}
            After computing the time derivative of  \eqref{eq:modified momentum NLSH},
            we use \eqref{eq:claim imaginary 1} and \eqref{eq:conjugated system} to obtain
            \begin{align}
                \frac{d}{dt} \bar{P}_j[Q\left(t\right)] &= 2\mathrm{Im} \int_{\Real^N} \frac{\partial}{\partial t} q_0^*\frac{\partial}{\partial x_j} q_0 +
                \tau \sum_{k=1}^{N} \frac{\partial}{\partial t} q_k^* \frac{\partial}{\partial x_j} q_k \, d\bfx\\
                & = 2\mathrm{Im} \int_{\Real^N} - i \mathrm{div}\left(\mathbf{q}_1^{*}\right)\frac{\partial}{\partial x_j} q_0 - i \kappa |q_0|^{p-1} q_0^{*}\frac{\partial}{\partial x_j} q_0  +
                i\sum_{k=1}^{N} \partial_k q_0^{*} \frac{\partial}{\partial x_j} q_k - i\sum_{k=1}^{N} q_k^{*}\frac{\partial}{\partial x_j} q_k  \, d\bfx \label{eq:dt modified momentum}
            \end{align}
            For the second term in \eqref{eq:dt modified momentum}, we have that
            \begin{align*}
                -2\kappa \mathrm{Im} \int_{\Real^N} i |q_0|^{p-1} q_0^{*}\frac{\partial}{\partial x_j} q_0 \, d\bfx 
                = -2\kappa  \int_{\Real^N} |q_0|^{p-1}\mathrm{Re} \left(q_0^{*}\frac{\partial}{\partial x_j} q_0\right) \, d\bfx = 
                -\frac{2\kappa}{p+1} \int_{\Real^N} \frac{\partial}{\partial x_j} |q_0|^{p+1} \, d\bfx
            \end{align*}
            Similarly, for the last term in \eqref{eq:dt modified momentum}, we get 
            \begin{align*}
                -2 \mathrm{Im} \int_{\Real^N} i \sum_{k=1}^{N} q_k^{*}\frac{\partial}{\partial x_j} q_k \, d\bfx =
                -2  \sum_{k=1}^{N} \int_{\Real^N} \mathrm{Re}\left(q_k^{*}\frac{\partial}{\partial x_j} q_k\right) \, d\bfx =
                    - \sum_{k=1}^{N} \int_{\Real^N} \frac{\partial}{\partial x_j} |q_k|^2 \, d\bfx.
            \end{align*}
            Finally, for the first and third terms in \eqref{eq:dt modified momentum}, we have 
            \begin{align*}
                2\mathrm{Im} \int_{\Real^N} - i \mathrm{div}\left(\mathbf{q}_1^{*}\right)\frac{\partial}{\partial x_j} q_0 + i\sum_{k=1}^{N} \partial_k q_0^{*} \frac{\partial}{\partial x_j} q_k \, d\bfx
                &= 2\sum_{k=1}^{N} \mathrm{Re} \left[\int_{\Real^N} q_k^{*}\frac{\partial^2}{\partial x_j \partial x_k} q_0 - \frac{\partial^2}{\partial x_k \partial x_j} q_0^{*} q_k \, d\bfx\right]. 
            \end{align*}
        \item For the conservation of the modified energy, it suffices to observe that
        \begin{align}
            \frac{d}{dt} \Hmod = \int_{\Real^N} \frac{\delta \Hmod}{\delta Q} \cdot Q_t + \frac{\delta \Hmod}{\delta Q^*} \cdot Q^*_t \, d\bfx
             = \int_{\Real^N} \left(\frac{\delta \Hmod}{\delta Q^*}\right)^* \cdot Q_t + \frac{\delta \Hmod}{\delta Q^*} \cdot Q^*_t \, d\bfx
             = 2\mathrm{Re} \int_{\Real^N} \left(\frac{\delta \Hmod}{\delta Q^*}\right)^* \cdot Q_t \, d\bfx,
        \end{align}
        since $\bar H$ is real-valued.
        Using the Hamiltonian structure of the NLSH system \eqref{eq:NLSH as a Hamiltonian}, we have that
        \begin{align}
            \frac{d}{dt} \Hmod = 2\mathrm{Re} \int_{\Real^N} \left(\frac{\delta \Hmod}{\delta Q^*}\right)^* \cdot \left(-i \bar{\mathcal{J}} \frac{\delta \bar{H}}{\delta Q^*}\right) \, d\bfx =
             2\mathrm{Re}\left[i \int_{\Real^N} \left|\frac{\delta \Hmod}{\delta q_0^*}\right|^2 + \tau^{-1}\left|\frac{\delta \Hmod}{\delta \bfqone^*}\right|^2\, d\bfx \right]  = 0.
        \end{align}
    \end{enumerate}
\end{proof}
\begin{remark}
    \label{remark: potential}
    For $V_0(\bfx) \neq 0$, it is straighforward to show that the modified mass $\bar{M}$ is still conserved by the NLSH system.
    The modified energy has to be redefined as $\Hmod_{V_0} = \Hmod - \int_{\Real^N} V_0(\bfx)|q_0|^2 d\bfx$ in order to be conserved,
    with the additional requirement that $V_0|q_0|^2 \in L^1(\Real^N)$.
    Momentum is no longer conserved, and we have instead the balance law
    \begin{align}
        \frac{d}{dt} \bar{P}_j[Q(t)] = -\int_{\Real^N} |q_0|^2 \frac{\partial}{\partial x_j} V_0(\bfx) \, d\bfx,
    \end{align}
    analogous to the NLS case.
\end{remark}

\section{Special solutions}
\label{sec:special}
\subsection{Plane waves and dispersion relation}
The NLS equation \eqref{eq:ND-NLS} 
admits plane wave solutions of the form
\begin{align}
    \label{eq:ND-plane-wave-NLS}
    u(\bfx,t) = \sigma e^{i (\bfk \cdot \bfx - \omega t)},
\end{align}
where $\sigma \in \Real$, $\mathbf{k} \in \Real^N$, and $\omega \in \Real$ satisfy the dispersion relation
\begin{align}
    \label{eq:ND-LDR-NLS}
    \omega(\bfk) = \ktk - \kappa \sigma^{p-1}.
\end{align}
Similarly, the NLSH system \eqref{eq:ND-NLSH} admits plane wave solutions of the form
\begin{align}
    \label{eq:ND-plane-wave-NLSH}
    q_0^\pm(\bfx,t) = \sigma e^{i (\mathbf{k} \cdot \bfx - \omega_\tau^\pm t)},\quad
    \mathbf{q}_1^\pm(\bfx,t) = \frac{i \bfk}{1 + \tau \omega_\tau^\pm} q_0^\pm(\bfx,t),
\end{align}
where $\omega_\tau^\pm \in \Complex$ satisfies the dispersion relation
\begin{align}
    \label{eq:ND-DR-NLSH}
    (\omega_\tau^\pm+\kappa\sigma^{p-1})(1+\tau \omega_\tau^\pm) = \ktk
\quad \Rightarrow \quad
    \omega^\pm_\tau(\bfk) = \frac{-1-\kappa\sigma^{p-1}\tau\pm\sqrt{4\tau \ktk  +(\kappa\sigma^{p-1}\tau-1)^2}}{2\tau}.
\end{align}

In the linear case ($\kappa = 0$), the dispersion relation \eqref{eq:ND-DR-NLSH} simplifies to
$
    \frac{-1\pm\sqrt{1+4\tau \ktk}}{2\tau}.
$
Thus, like for NLS, the trivial solution $\bfq=\mathbf{0}$ of the NLSH system is linearly stable since $\mathrm{Im}(\omega_\tau^\pm) = 0$ for all $\bfk \in \Real^N$.
Furthermore, we have that
 \begin{align*}
    \omega_\tau^{+}(\bfk) = \omega(\bfk)-\tau(\ktk-\kappa\sigma^{p-1})^2+\bigo(\tau^2),
 \end{align*}
hence   $\omega_\tau^{+}(\bfk) \to \omega(\bfk)$ as $\tau\to 0$.
On the other hand, $\omega_\tau^{-}$, which corresponds to a spurious mode associated with the hyperbolic relaxation,
diverges like $-1/\tau$ as $\tau \to 0$.

\subsubsection{Modulational instability of NLSH}
\label{sec:MI theoretical}
Let us now study the stability of the plane waves \eqref{eq:ND-plane-wave-NLSH} for $\bfk=\bfzero$.
Then, the dispersion relation \eqref{eq:ND-DR-NLSH} simplifies to
\begin{align}
    \label{eq:ND-DR-NLSH-k0}
    \omega^+_\tau(\bfzero) = -\kappa \sigma^{p-1}, \quad
    \omega^-_\tau(\bfzero) = -\frac{1}{\tau},
\end{align}
and the corresponding physical plane wave solution (omitting the superscript $+$) is given by
\begin{align}
    \label{eq:ND-plane-wave-NLSH-k0}
    q_0(\bfx,t) = \sigma e^{i \kappa \sigma^{p-1} t}, \quad
    \mathbf{q}_1(\bfx,t) = \mathbf{0}.
\end{align}
We consider a perturbation of the form
\begin{align}
    \label{eq:ND-plane-wave-NLSH-k0-perturbation}
    q_0(\bfx,t) = \left(\sigma + \eps(\bfx,t)\right) e^{i\kappa \sigma^{p-1} t}, \quad 
    \bfqone(\bfx,t) = e^{i\kappa \sigma^{p-1} t} \bfdelta(\bfx,t),
\end{align}
where $\eps(\bfx,t)\ll \sigma$ and $|\bfdelta(\bfx,t)|\ll 1$ are allowed to be complex-valued.
Substituting \eqref{eq:ND-plane-wave-NLSH-k0-perturbation} into \eqref{eq:ND-NLSH} and using the linearization
\begin{align*}
|\sigma+\eps|^{p-1}(\sigma+\eps) = \sigma^p + \tfrac{p+1}{2}\sigma^{p-1}\eps + \tfrac{p-1}{2}\sigma^{p-1}\eps^* + O(\eps^2),
\end{align*}
we get the following linear system for the perturbation:
\begin{subequations}
\begin{align}
    \label{eq:ND-NLSH-k0-perturbation-eq1}
    i\eps_t + \mathrm{div}(\bfdelta) + \frac{p-1}{2}\kappa \sigma^{p-1} \left(\eps + \eps^*\right) = 0,\\
    \label{eq:ND-NLSH-k0-perturbation-eq2}
    i \tau \bfdelta_t - \nabla \eps + (1-\tau\kappa \sigma^{p-1})\bfdelta = 0.
\end{align}
\end{subequations}
Then, we consider a sideband ansatz of the form
\begin{align}
    \label{eq:sideband ansatz NLSH}
    \eps(\bfx,t) = 
    \eps_+ e^{i(\bfk\cdot \bfx - \Omega t)} + \eps_-^* e^{-i(\bfk\cdot \bfx - \Omega^* t)},\quad 
    \bfdelta(\bfx,t) = \bfdelta_+ e^{i(\bfk\cdot \bfx - \Omega t)} + \bfdelta_-^* e^{-i(\bfk\cdot \bfx - \Omega^* t)},
\end{align}
where $\eps_+$, $\eps_-$, $\bfdelta_+$, and $\bfdelta_-$ are complex constants, $\bfk \in \Real^N$ is the wavevector of the perturbation,
and $\Omega \in \Complex$ is the frequency of the perturbation.
After substituting \eqref{eq:sideband ansatz NLSH} into \eqref{eq:ND-NLSH-k0-perturbation-eq2}, we 
find the following algebraic relation between $\bfdelta_\pm$ and $\eps_\pm$:
\begin{align}
    \label{eq:algebraic relation delta-eps}
    \bfdelta_\pm = \frac{i \bfk \eps_\pm}{1-\tau\kappa\sigma^{p-1} \pm \tau \Omega}.
\end{align}
Then, after noting that
\begin{subequations}
\begin{gather}
i\partial_t \eps = \Omega \eps_+ e^{i(\bfk\cdot \bfx - \Omega t)} -  \Omega^*\,\eps_-^* e^{-i(\bfk\cdot \bfx - \Omega^* t)},
\quad
\mathrm{div}(\bfdelta) = i\bfk^T \left(\boldsymbol\delta_+e^{i(\bfk\cdot x-\Omega t)}
-\boldsymbol\delta_-^*e^{-i(\bfk\cdot x-\Omega^*t)}\right)
\\
(\eps+\eps^*) = (\eps_+ + \eps_-)e^{i(\bfk\cdot \bfx - \Omega t)} + (\eps_+^* + \eps_-^*) e^{-i(\bfk\cdot \bfx - \Omega^* t)},
\end{gather}
\label{eq:intermediate perturbation linear system}
\end{subequations}
and using \eqref{eq:algebraic relation delta-eps}, 
we can collect terms proportional to $e^{i(\bfk\cdot \bfx - \Omega t)}$ and $e^{-i(\bfk\cdot \bfx - \Omega^* t)}$ in
 \eqref{eq:ND-NLSH-k0-perturbation-eq1} and its complex conjugate to get the following linear system for $\eps_+$ and $\eps_-$:
\begin{align*}
\begin{bmatrix}
    \Omega - \frac{\ktk}{1-\tau\kappa\sigma^{p-1} + \tau \Omega} + \frac{p-1}{2}\kappa \sigma^{p-1} & \frac{p-1}{2}\kappa \sigma^{p-1} \\
    \frac{p-1}{2}\kappa \sigma^{p-1} & -\Omega - \frac{\ktk}{1-\tau\kappa\sigma^{p-1} - \tau \Omega} + \frac{p-1}{2}\kappa \sigma^{p-1} 
\end{bmatrix}
\begin{bmatrix}
    \eps_+ \\
    \eps_-
\end{bmatrix} = 0.
\end{align*}
A nontrivial solution exists if the determinant of the matrix is zero.
To ease the algebra, let's introduce the following notation:
\begin{align*}
    a = 1-\tau\kappa\sigma^{p-1}, \quad
    b = \frac{p-1}{2}\kappa \sigma^{p-1}, \quad
    r = a+\tau \Omega, \quad
    s = a-\tau \Omega.
\end{align*}
Then, the determinant condition can be written as
\begin{gather}
    \label{eq:ND-NLSH-k0-perturbation-determinant}
    rs(b^2-\Omega^2)-r(\Omega+b)\ktk-s(b-\Omega)\ktk+|\bfk|^4 -rsb^2 = 0.
\end{gather}
Noting that $rs = a^2 - \tau^2 \Omega^2$, $r-s = 2\tau \Omega$, and $r+s = 2a$, Equation \eqref{eq:ND-NLSH-k0-perturbation-determinant} 
becomes a quadratic equation in $\Omega^2$:
\begin{align}
    \label{eq:ND-NLSH-k0-perturbation-quartic}
    \tau^2 \Omega^4 - \left(a^2+2\tau\ktk\right)\Omega^2 + \ktk\left(\ktk - 2ab\right) = 0.
\end{align}
The roots of this equation are given by
\begin{align*}
    \Omega_\pm^2 = \frac{a^2 + 2\tau \ktk \pm \sqrt{D}}{2\tau^2},
\end{align*}
where 
\begin{align}
    \label{eq:ND-NLSH-k0-perturbation-discriminant}
    D 
    &= a^4 + 4\tau \ktk a^2 + 8 \tau^2 ab \ktk.
\end{align}
The plane wave solution \eqref{eq:ND-plane-wave-NLSH-k0} is modulationally unstable if $\mathrm{Im}(\Omega_\pm) > 0$.
This can occur through two different mechanisms: 
\begin{itemize}
    \item ($D>0$ and $\Omega_\pm^2 < 0$) In this case, we have that
        $
        \Omega^2_+ +\Omega^2_-= \tau^{-2}({a^2 + 2\tau \ktk}) > 0,
        $
        so one of the roots is negative if $\Omega_-^2 \Omega_+^2  = {\tau^{-2}} {\ktk\left(\ktk - 2ab\right)}< 0$, i.e.,
        \begin{align}
            \label{eq: MI condition NLSH}
            \ktk < 2ab = (p-1)\kappa \sigma^{p-1}(1-\tau\kappa\sigma^{p-1}).
        \end{align}
        This reduces to the classical Benjamin-Feir condition for NLS in the limit $\tau \to 0$.
        Moreover, in the focusing case ($\kappa > 0$), condition \eqref{eq: MI condition NLSH} gives us a band of unstable wavevectors
        whose width decreases like $\sqrt{1-\tau\kappa\sigma^{p-1}}$.
        In the defocusing case ($\kappa < 0$), analogously to NLS, \eqref{eq: MI condition NLSH}  is never satisfied.
    \item ($D < 0$) It is easy to see that this condition arises uniquely from the hyperbolic relaxation ($\tau>0$),
        leading to new band of unstable wavevectors, given by
        \begin{align}
            \label{eq: MI condition NLSH D<0}
            \ktk > \frac{-a^3}{4\tau a+ 8 \tau^2 b}\eqqcolon \theta(\kappa\sigma^{p-1}, \tau),\quad  \text{and}  \quad 4\tau a^2 + 8 \tau^2 ab < 0.
        \end{align}
        Condition \eqref{eq: MI condition NLSH D<0} gives us a practical restriction
        on the choice of the relaxation parameter $\tau$ for which NLSH
        faithfully reproduces the modulational stability properties of NLS.
        In the focusing case ($\kappa > 0$),
        the instability is suppressed for every $\ktk$ if $\tau < \kappa^{-1}\sigma^{1-p}$,
        while in the defocusing case ($\kappa < 0$), it suffices to take $\tau<\kappa^{-1}\sigma^{1-p}/(2-p)$ to avoid the instability.
        This is further motivated and verified numerically in Section \ref{sec:MI numerical}.
        The lower bound on the unstable band \eqref{eq: MI condition NLSH D<0} is illustrated
        in Figure \ref{fig:MI bands} for the focusing and defocusing cases for different values of $\tau$
        and $\kappa \sigma^{p-1}$.
        \begin{figure}
            \centering
            \includegraphics[width=0.8\textwidth]{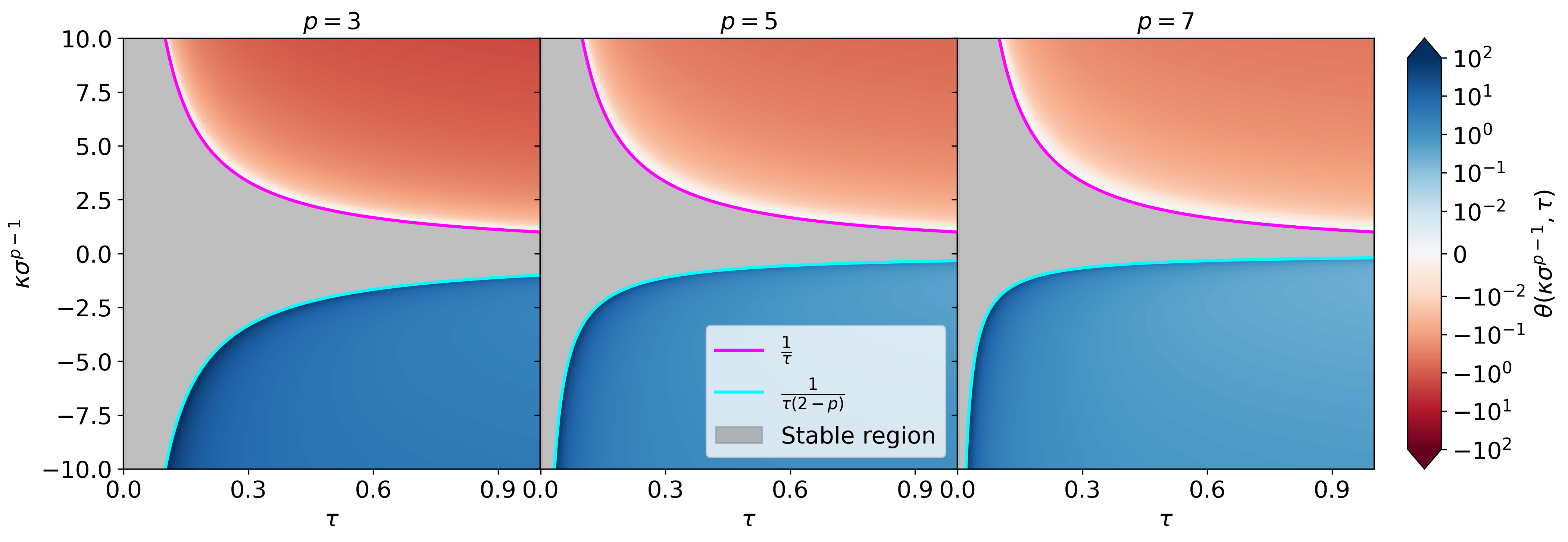}
            \caption{Lower bound of the hyperbolization-induced unstable band for NLSH $\theta$ (see Equation \eqref{eq: MI condition NLSH D<0}) in the focusing and defocusing cases for different
            values of $\tau$ and $\kappa \sigma^{p-1}$.
            The instability is suppressed for every $\ktk$ in the gray region, while it is unavoidable for every $\ktk$ in the red region,
            where $\theta\leq0$.
            }
            \label{fig:MI bands}
        \end{figure}
\end{itemize}
Note that the analysis in this section is restricted to plane waves with $\bfk=\bfzero$.
In the case of NLS, this extends directly to arbitrary $\bfk$ due to the Galilean invariance
of the system.
Such a Galilean boost is not available for NLSH (although
the system supports non-trivial traveling wave solutions, to be reported in a forthcoming paper),
thus a different argument would be required in the general case, which 
is left for future work.

\subsection{Standing wave solutions}
Let $\kappa>0$, $\mu >0$, and $V_0(\bfx) = 0$.
Standing wave solutions of the NLS equation \eqref{eq:ND-NLS} are given by
\begin{align}
    u(\bfx,t) = e^{i \mu t} \tilde{u}(\bfx),
\end{align}
where $\tilde{u}$ satisfies the elliptic equation
\begin{align}
    \label{eq:steady-NLS}
    \nabla^2 \tilde{u} - \mu \tilde{u} + \kappa |\tilde{u}|^{p-1} \tilde{u} = 0.
\end{align}
Solutions of \eqref{eq:steady-NLS} are well-studied in the literature.
In particular, if $\tilde{u}\neq 0$ and $1<p<\frac{N+2}{N-2}\; (p>1$ if $N=1,2)$,
then there exists a unique (up to translations)
positive, radially symmetric, decaying solution of \eqref{eq:steady-NLS},
$\overline{u}\in H^1(\Real^N)$,
which minimizes the energy functional \eqref{eq:energy NLS}
over all its critical points, see e.g. \cite[Thms. 8.1.4 to 8.1.6]{cazenave2003semilinear}.
$\overline{u}$ is referred to as the ground state solution of \eqref{eq:steady-NLS}.

Similarly, for some $\tau>0$, we look for standing wave solutions of the NLSH system \eqref{eq:ND-NLSH} 
of the form
\begin{align*}
 \begin{pmatrix}
    q_0\\
    \bfqone
 \end{pmatrix}(\bfx,t)
 =
    e^{i \mu t} \begin{pmatrix}
        \qtzero\\
        \bfqtone
    \end{pmatrix}(\bfx),
\end{align*}
Substituting this ansatz into \eqref{eq:ND-NLSH}, we obtain the system
\begin{subequations}
\begin{align}
    -\mu \qtzero + \mathrm{div}(\bfqtone) + \kappa {|\qtzero|}^{p-1}\qtzero  = 0, \label{eq:steady-NLSH 1}\\
    (1-\tau \mu) \bfqtone  = \nabla \qtzero . \label{eq:steady-NLSH 2}
\end{align}
\end{subequations}
Replacing \eqref{eq:steady-NLSH 2} into \eqref{eq:steady-NLSH 1}, we obtain
\begin{align}
    \label{eq:steady-NLSH}
    \frac{1}{1-\tau \mu} \nabla^2 \qtzero  -\mu\qtzero +\kappa {|\qtzero|}^{p-1}\qtzero= 0.
\end{align}
It is immediate to see that, if $1-\tau \mu>0$, then there exists a unique (up to translations)
positive, radially symmetric, decaying solution of \eqref{eq:steady-NLSH},
$\overline{q}_0\in H^1(\Real^N)$, given by
\begin{align*}
    \overline{q}_0(\bfx) = \overline{u}\left(\sqrt{1-\tau \mu} \bfx \right),
\end{align*}
which minimizes the modified energy functional \eqref{eq:modified energy NLSH}
over all its critical points.
Together with \eqref{eq: MI condition NLSH D<0}, $1-\tau\mu>0$
is an additional smallness condition on $\tau$
for a faithful approximation of the NLS dynamics with NLSH.
Moreover, the convergence result presented in one dimension in \cite[Prop. 1]{biswas2025hyperbolic}
can be extended to $N$ dimensions as follows.
\begin{proposition}
    \label{prop: converge ground states}
    Let $\kappa>0$, $\mu>0$, and $1<p<\frac{N+2}{N-2}\; (p>1$ if $N=1,2)$.
    Then the family \eqref{eq:steady-NLSH} of ground state solutions $\{\overline{q}_0\}_{\tau>0}$ of the NLSH system converges uniformly,
    and linearly in $\tau$,
    to the ground state solution $\overline{u}$ of the NLS equation \eqref{eq:steady-NLS}  as $\tau \to 0$.
    Furthermore, $\{\overline{\bfq}_1\}_{\tau>0}=\{(1-\tau \mu )^{-1}\nabla\overline{q}_0\}_{\tau>0}$ converges uniformly and linearly in $\tau$ to $\nabla \overline{u}$.
\end{proposition}
\begin{proof}
    The result follows from the explicit relation between $\overline{q}_0$  and $\overline{u}$ \eqref{eq:steady-NLSH},
    the use of the Mean Value Theorem, and the fact that $\nabla \overline{u}$ decays exponentially at infinity
    (see e.g. \cite[Thm. 8.1.1]{cazenave2003semilinear}).
    The convergence of $\{\overline{\bfq}_1\}_{\tau>0}$  to $\nabla \overline{u}$ can be obtained following a similar argument,
    using the exponential decay of higher order derivatives of $\overline{u}$.
\end{proof}

\section{Convergence of hyperbolized solutions to the NLS solution} \label{sec:convergence}
We prove the convergence of the hyperbolic approximation of the NLS system to the NLS equation for smooth solutions using the relative energy stability framework,
following \cite{giesselmann2025convergence}. The key idea is to view a smooth solution of the NLS equation as satisfying a perturbed form of the NLSH system, where 
the perturbation terms are proportional to $\tau$. Then we use the relative energy stability framework for the NLSH system to estimate the difference between the 
exact solution of the NLSH system and the solution of the NLS equation in terms of these perturbation terms, which allows us to conclude the convergence of the 
hyperbolized solution to that of the original NLS equation.
\begin{theorem}
    \label{thm: convergence ND}
Let $u\in H^3((0,T)\times\mathbb R^N)\cap L^\infty((0,T)\times\mathbb R^N)$ be a solution to \eqref{eq:ND-NLS} with initial data $u_0\in H^3(\mathbb R^N)$. 
For each $0<\tau\le\tau_0$, where $\tau_0>0$ is fixed, let $q\in H^1((0,T)\times\mathbb R^N)$ be a solution to \eqref{eq:ND-NLSH} with initial data $q|_{t=0}=[u_0,\nabla u_0]^T$. 
Assume that
\begin{equation}
\sup_{0<\tau\le\tau_0}
\|q_0\|_{L^\infty((0,T)\times\mathbb R^N)}<\infty.
\end{equation}
Then
\begin{align}
    \|q_0-u\|_{L^\infty([0,T]; L^2(\mathbb{R}^N))} + \sum \limits_{j=1}^{N}\sqrt{\tau} \, \|q_j - \partial_{x_j} u\|_{L^\infty([0,T]; L^2(\mathbb{R}^N))} = \mathcal{O}(\tau) \;.
\end{align}
\end{theorem}
\begin{proof}
We consider the NLSH system \eqref{eq:ND-NLSH} with a real-valued potential $V_0(x)$:
\begin{subequations}
    \label{eq:ND-NLSH_1}
    \begin{align}
      i \partial_t q_0 + \nabla \cdot \mathbf{q}_1 + (V_0(x)+\kappa |q_0|^2)q_0 = 0,\\
   i \tau \partial_t \mathbf{q}_1 = \nabla q_0 - \mathbf{q}_1.
    \end{align}  
\end{subequations}
We define an approximate solution to the hyperbolized system \eqref{eq:ND-NLSH_1} as
\begin{subequations} \label{Eq:ND-NLSH_approx_sol}
\begin{align}
\qt_0 &= u\;, \\
\bqt_1 &= \nabla u - i \tau \, \partial_t (\nabla u) \;.
\end{align}
\end{subequations}
Substituting the approximate solution $\qt = [\qt_0,\bqt_1 ]^T$ into \eqref{eq:ND-NLSH_1} shows that $\qt$ satisfies a perturbed NLSH system given by 
\begin{subequations} \label{Eq:Pert_ND-NLSH}
\begin{align}
    i\partial_t \qt_0 + \nabla \cdot \bqt_1 + (V_0(x)+\kappa |\qt_0|^2)\qt_0 &=  -i \tau R_0\;, \label{Eq:pert_NLSH1} \\
    i \tau \partial_t \bqt_1 &= \nabla \qt_0 - \bqt_1 + \tau \bR_1 \;, \label{Eq:pert_NLSH2}
\end{align}
\end{subequations}
where the residuals are
\begin{align*}
R_0 &= \partial_t(\Delta u)
     = \partial_t\left(
       \sum_{j=1}^{N}\frac{\partial^2 u}{\partial x_j^2}
       \right),\\
\bR_1 &= \tau\partial_{tt}(\nabla u)
      = \tau(\partial_{ttx_1}u,\ldots,
             \partial_{ttx_N}u)^T.
\end{align*}
Consider the relative energy 
\begin{align}\label{Eq:ND-rel_energy}
    \eta(q, \qt) = \frac{1}{2} \int_{\Real^N} \left( |\delta q_0|^2 + \tau |\delta \bq_1|^2 \right) \, dx = \frac{1}{2} \int_{\Real^N} \left( |\delta q_0|^2 + \tau \sum\limits_{j=1}^N|\delta q_j|^2  \right) \, dx
\end{align}
where $\delta q_0 = q_0 - \qt_0$ and $\delta \bq_1 = \bq_1 - \bqt_1$.
The time evolution of the relative energy is given by
\begin{align}\label{Eq:ND-rel_energy_1}
    \frac{d}{dt} \eta(q,\qt) = \operatorname{Re} \int_{\Real^N} \left(\delta q_0^* \cdot \partial_t (\delta q_0) + \tau \sum \limits_{j=1}^{N}\delta q_j^* \cdot \partial_t (\delta q_j) \right) \, dx \;.
\end{align} 
Subtracting \eqref{Eq:Pert_ND-NLSH} from \eqref{eq:ND-NLSH_1} we obtain
\begin{align*}
    \partial_t (\delta q_0) & = i \left(\nabla \cdot (\delta \bq_1)+ V_0(x)\delta q_0 + \kappa (|q_0|^2 q_0 - |\qt_0|^2 \qt_0)\right) + \tau R_0 \;,\\
    \tau \partial_t (\delta \bq_1) & = -i \nabla (\delta q_0) + i \delta \bq_1 + i \tau \bR_1.
\end{align*}
Substituting these expressions into the relative energy time evolution equation \eqref{Eq:ND-rel_energy_1} and integrating by parts over $\mathbb R^N$, we find
\begin{align}\label{Eq:ND-rel_energy_2}
\begin{aligned}
    \frac{d}{dt} \eta(q,\qt)
    = \operatorname{Re} \int_{\mathbb{R}^N} \Big(
        i \delta q_0^{*}\, \nabla \cdot (\delta \bq_1)
        + i V_0(x)|\delta q_0|^2
        + i \kappa\, \delta q_0^{*} (|q_0|^2 q_0 - |\qt_0|^2 \qt_0)
        + \tau \delta q_0^{*}\,  R_0  \\
        {} \qquad\qquad\qquad
        +\, i \delta q_0\, \nabla \cdot (\delta \bq_1^{*})
        + i \sum_{j=1}^N |\delta q_j|^2
        + i \sum_{j=1}^N \tau\, \delta q_j^{*} R_j
    \Big)\, dx \;.
\end{aligned}
\end{align}
Observe that $i \delta q_0^* \nabla \cdot (\delta \bq_1)
+ i \delta q_0 \nabla \cdot (\delta \bq_1^*)
= 2i \operatorname{Re}
\left(\delta q_0^* \nabla \cdot (\delta \bq_1)\right)$
and $i\sum_{j=1}^N|\delta q_j|^2$ are purely imaginary.
Moreover, since $V_0(x)$ is real-valued,
$iV_0(x)|\delta q_0|^2$ is also purely imaginary. 
Thus, these terms have zero real part, and \eqref{Eq:ND-rel_energy_2} reduces to
\begin{align}
    \frac{d}{dt} \eta(q,\qt) & = \operatorname{Re} \int_{\Real^N} \left(i \kappa \delta q_0^* (|q_0|^2 q_0 - |\qt_0|^2 \qt_0) + \tau \delta q_0^*  R_0 
    + i \sum_{j=1}^N \tau\, \delta q_j^{*} R_j \right) \, dx \;, \\
     & = \operatorname{Re} \int_{\Real^N} i \kappa \delta q_0^* (|q_0|^2 q_0 - |\qt_0|^2 \qt_0) \, dx + \operatorname{Re} \int_{\Real^N} \tau \delta q_0^*  R_0 \, dx 
     + \sum \limits_{j=1}^N \operatorname{Re} \int_{\Real^N} i \tau \delta q_j^*  R_j  \, dx \;. \label{Eq:ND-rel_energy_3}
\end{align}
By the assumptions of the theorem, both $q_0$ and
$\qt_0=u$ are uniformly bounded. Define
\[
M=\max\left\{
\|u\|_{L^\infty((0,T)\times\mathbb R^N)},
\sup_{0<\tau\le\tau_0}
\|q_0\|_{L^\infty((0,T)\times\mathbb R^N)}
\right\}.
\]
Then $M$ is finite and independent of $\tau$, and
$|q_0|,|\qt_0|\le M$ almost everywhere.
Consequently,
$||q_0|^2q_0-|\qt_0|^2\qt_0|
\le 3M^2|\delta q_0|$, and hence
\begin{align}\label{Eq:ND-term_1}
    \operatorname{Re} \int_{\Real^N} i \kappa \delta q_0^* (|q_0|^2 q_0 - |\qt_0|^2 \qt_0) \, dx \leq 6 |\kappa| M^2 \eta(q,\qt) \;.
\end{align}
Using Young's inequality, we obtain
\begin{align}\label{Eq:ND-term_23}
    \operatorname{Re}\int_{\Real^N}
    \delta q_0^* \tau R_0\,dx
    &\le
    \frac{1}{2}\|\delta q_0\|_{L^2}^2
    +\frac{\tau^2}{2}
    \|\partial_t(\Delta u)\|_{L^2}^2,
    \\
    \operatorname{Re}\int_{\Real^N}
    i\delta q_j^*\tau R_j\,dx
    &\le
    \frac{\tau}{2}\|\delta q_j\|_{L^2}^2
    +\frac{\tau^3}{2}
    \|\partial_{ttx_j}u\|_{L^2}^2,
    \qquad j=1,\ldots,N.
\end{align}
Using the estimates \eqref{Eq:ND-term_1} and
\eqref{Eq:ND-term_23} in \eqref{Eq:ND-rel_energy_3},
we obtain
\begin{align}\label{Eq:ND-rel_energy_4}
    \frac{d}{dt}\eta(q,\qt)
    &\le
    (6|\kappa|M^2+1)\eta(q,\qt)
    +\frac{\tau^2}{2}
    \|\partial_t(\Delta u)\|_{L^2}^2
    +\frac{\tau^3}{2}
    \sum_{j=1}^{N}
    \|\partial_{ttx_j}u\|_{L^2}^2 \;.
\end{align}
By the initial conditions in the theorem,
$q|_{t=0}=[u_0,\nabla u_0]^T$.
Since $\qt_0|_{t=0}=u_0$ and
$\bqt_1|_{t=0}=\nabla u_0
-i\tau\partial_t(\nabla u)|_{t=0}$,
we obtain
$\eta(q,\qt)|_{t=0}=\mathcal O(\tau^3)$.
Applying Gronwall's lemma to
\eqref{Eq:ND-rel_energy_4}, we conclude that
\begin{align}
    \|\eta(q, \qt)\|_{L^\infty([0,T])} \le C_T(\tau^2+\tau^3) =\mathcal O(\tau^2)\;,
\end{align}
where $C_T>0$ is independent of $\tau$. 
Since $\widetilde q_0=u$ and
$\widetilde{\mathbf q}_1=\nabla u-i\tau\nabla u_t$,
the preceding estimate implies
\[
\|q_0-u\|_{L^\infty(0,T;L^2(\mathbb R^N))}
=\mathcal O(\tau).
\]
Moreover, by the triangle inequality,
\[
\begin{aligned}
\sqrt{\tau}\,
\|\mathbf q_1-\nabla u\|_{L^\infty(0,T;L^2(\mathbb R^N))} \le \sqrt{\tau}\, \|\mathbf q_1-\widetilde{\mathbf q}_1\|_{L^\infty(0,T;L^2(\mathbb R^N))} +
\tau^{3/2}\|\nabla u_t\|_{L^\infty(0,T;L^2(\mathbb R^N))} =\mathcal O(\tau).
\end{aligned}
\]
\end{proof}

\begin{remark}
    The foregoing result is of both practical and theoretical significance.  Nevertheless, the
    assumptions it uses cannot be taken for granted in general.  Numerical experiments strongly
    suggest that, for large enough $\tau$ and large enough initial data, solutions of NLSH can
    form shock discontinuities.  The same experiments suggest (reassuringly) that such discontinuities
    are avoided when $\tau$ is small enough.  Furthermore, it is well known that some solutions of NLS
    blow up in $L^\infty$ (see e.g. \cite{yang2010nonlinear}); the theorem cannot be applied arbitrarily close to the blowup time in such cases.
\end{remark}

\subsection{Well-prepared initial data}
    \label{sec: well-prepared initial data}
    For $\tau > 0$, we define the slow manifold of the NLSH system as
    \begin{align}
        \slowmani = \left\{ (q_0, \bfqone) \in H^1(\mathbb{R}^N\times \Real,
        \mathbb{C})\times L^2(\mathbb{R}^N\times \Real, \mathbb{C}^N) : \bfqone(\cdot, t) =
         \slowmap(\cdot,t) \; \text{such that \eqref{eq:ND-NLSH} holds} \right\},
    \end{align}
    i.e., the set of all the solutions to NLSH $\bfq=[q_0, \bfqone]^\top$ such that $\bfqone$ is slaved to $q_0$ through the map $\mathcal{M}_\tau$.
    For $\bfq \in \slowmani$, we define $\induced := i(\mathrm{div} \slowmap + V_0 q_0+ \kappa |q_0|^2 q_0) = \partial_t q_0$.
    Then, using \eqref{eq:ND-NLSH} and the chain rule,  we can characterize the slow manifold through the following fixed-point problem
    \begin{align}
        \label{eq: fixed-point problem}
        \slowmap = \nabla q_0 - i \tau \braket{D\slowmap,\induced}
    \end{align}
    where $D\slowmap$ is the Fréchet derivative of $\mathcal{M}_\tau$ at $q_0$.

    If we make the (strong) assumption that $\braket{D\slowmap,\induced}$ is uniformly bounded in $\tau$,
    we get the $\bigo(\tau)$ approximation of $\slowmap$, $\smapiter{0}=\nabla q_0$.
    We can think of Theorem \ref{thm: convergence ND} as proof that the solution to NLSH
    with initial data $\bfqone(\cdot,0)=\smapiter{0}(\cdot,0)=\nabla u_0$ is $\bigo(\tau)$-close
    to the slow manifold of the NLSH system \eqref{eq:ND-NLSH} and, therefore, $\bigo(\tau)$-close
    to the solution of the NLS equation \eqref{eq:ND-NLS}.
    This invites the question of whether we can find a better approximation of the map $\slowmap$ to initialize the NLSH system.
    A natural approach is to apply fixed-point iterations to \eqref{eq: fixed-point problem} to obtain a sequence of approximations $\smapiter{n}$ of the slow manifold $\slowmap$,
    using $\smapiter{0}=\nabla q_0$ as the initial guess.
    Since $\mathcal{M}^0_\tau$ is linear, $\braket{D\smapiter{0},h} = \nabla h$ and we get
    \begin{align}
        \smapiter{1} = \nabla q_0 - i \tau \braket{D\smapiter{0},\inditer{0}}=
        \nabla q_0 + \tau \nabla (\nabla ^2 q_0 + V_0 q_0 + \kappa |q_0|^2 q_0).
    \end{align}
    For the next iteration, we have
    \begin{align}
        \smapiter{2} = \nabla q_0 - i \tau \braket{D\smapiter{1},\inditer{1}},
    \end{align}
    where $\inditer{1} = i(\mathrm{div} \smapiter{1} + V_0+ \kappa |q_0|^2 q_0)$, and
    \begin{align}
        \braket{D\smapiter{1},h} = \nabla h - i \tau \nabla (\nabla^2 h + \braket{DF[q_0],h}).
    \end{align}
    Here,  $F[q_0] = V_0 q_0 + \kappa |q_0|^2 q_0$ and
    \begin{align}
        DF[q_0]h = V_0 h + \kappa \frac{p+1}{2}|q_0|^{p-1} h + \kappa \frac{p-1}{2}|q_0|^{p-3} q_0^2 h^*.
    \end{align}

    Having computed an approximation $\smapiter{n}$ of the mapping $\slowmap$, 
    we can initialize the NLSH system with $n$-th order \emph{well-prepared initial data} $\bfqone(\cdot,0)=\smapiter{n}(\cdot,0)$.
    Moreover, substituting $\smapiter{n}$ into \eqref{eq:ND-NLSH} gives us a modified NLS equation for $q_0$
    that NLSH approximates to $\bigo(\tau^{n+1})$ accuracy if the initial data is well-prepared to order $n$.
    For instance, using $\smapiter{0}$ we get NLS, and using $\smapiter{1}$ we get the modified NLS equation
    \begin{align}
        i \partial_t q_0 +  (\mathcal{I}+ \tau \nabla^2) (\nabla^2 q_0 + V_0 q_0 + \kappa |q_0|^2 q_0) = 0. 
    \end{align}

\section{Numerical Methods} \label{sec:FVM}

In this and the next section, we leverage the extensive mathematical and computational
tools that have been developed for first-order hyperbolic PDEs in order to discretize
the NLSH system \eqref{eq:ND-NLSH} using a finite volume method.  This has several advantages.
Problems with sharp gradients, such as those arising in the study of dispersive shock waves and 
discontinuous trapping potentials can be handled robustly.
Moreover, well-established techniques like patch-based adaptive mesh refinement (AMR) and mapped grids 
can be used to accurately model a wide range of spatio-temporal scales and complex geometries.

For this, it is convenient to split the NLSH system into real and imaginary parts.
Letting $q_0 = z_0+i w_0$ and $\bfqone= \bfzone + i \bfwone$, we can rewrite the NLSH system as
\begin{align}
    \label{eq:ND-NLSH-real-imaginary}
\partial_t\left[\begin{array}{c}
z_0 \\
\bfzone \\
w_0 \\
\bfwone
\end{array}\right]+
\left[\begin{array}{c}
\operatorname{div}\left(\bfwone\right) \\
-\tau^{-1} \nabla w_0 \\
-\operatorname{div}\left(\bfzone\right) \\
\tau^{-1} \nabla z_0
\end{array}\right]=
\left[\begin{array}{l}
-\kappa\left(z_0^2+w_0^2\right) w_0-V_0(\bfx) w_0 \\
-\tau^{-1} \bfwone \\
\kappa\left(z_0^2+w_0^2\right) z_0+V_0(\bfx) z_0 \\
\tau^{-1} \bfzone
\end{array}\right].
\end{align}
Denoting by $\eta=\left[
z_0,
\bfzone,
w_0,
\bfwone \right]^T
$
the vector of real-valued unknowns, 
the flux terms can be expressed as
\begin{align}
\sum_{j=1}^{N} B_j \partial_{x_j} \eta  = 
\sum_{j=1}^{N}
 \begin{bmatrix}
0 & \mathbf{0}_{1\times N} & 0 & \mathbf{e}_{j}^T \\
\mathbf{0}_{N\times 1} & \mathbf{0}_{N\times N}  & -\tau^{-1}\mathbf{e}_{j} & \mathbf{0}_{N\times N}  \\
0 & -\mathbf{e}_{j}^T & 0 & \mathbf{0}_{1\times N} \\
\tau^{-1}\mathbf{e}_{j}  & \mathbf{0}_{N\times N}  & \mathbf{0}_{N\times 1} & \mathbf{0}_{N\times N} 
\end{bmatrix}
\partial_{x_j} \eta.
\end{align}
Note that we can write $B_j = P_j \tilde{B}_j P_j^T$,
where the projection matrix $P_j=[\mathbf{e}_{1},\mathbf{e}_{j+1},\mathbf{e}_{N+2},\mathbf{e}_{N+2+j}]$ 
has size $(2N+2)\times 4$ and $\tilde{B}_j$ is given by
\begin{align}
\tilde{B}_j = \begin{bmatrix}
0 & 0 & 0 & 1 \\
0 & 0 & -\tau^{-1} & 0 \\
0 & -1 & 0 & 0 \\
\tau^{-1} & 0 & 0 & 0
\end{bmatrix}.
\end{align}
Therefore, the nonzero eigenvalues of $B_j$ are equal (up to multiplicities) to the eigenvalues of $\tilde{B}_j$: $\lambda_{1,2} = -\tau^{-1/2}$,
$\lambda_{3,4} = \tau^{-1/2}$. Furthermore, the corresponding right eigenvectors  are given by 
$V=[v_1,v_2,v_3,v_4]=P_j\tilde{V}$ where $\tilde{V}$ is the matrix of right eigenvectors of $\tilde{B}_j$.
Explicitly, we have
\begin{align}
v_1 = \begin{bmatrix}
    -\sqrt{\tau}\\
    \mathbf{0}_{N\times 1}\\
    0\\
    \mathbf{e}_j
\end{bmatrix}, \quad
v_2 = \begin{bmatrix}
    0\\
    \sqrt{\tau^{-1}} \mathbf{e}_j\\
    1\\
    \mathbf{0}_{N\times 1}
\end{bmatrix}, \quad
v_3 = \begin{bmatrix}
    \sqrt{\tau}\\
    \mathbf{0}_{N\times 1}\\
    0\\
    \mathbf{e}_j
\end{bmatrix}, \quad
v_4 = \begin{bmatrix}
    0\\
    -\sqrt{\tau^{-1}} \mathbf{e}_j\\
    1\\
    \mathbf{0}_{N\times 1}
\end{bmatrix}.
\end{align}
The reader might recognize that, in absence of a source term, the hyperbolic part of the system becomes
two decoupled acoustics systems, one for $(z_0, \bfwone)$ and one for $(w_0, \bfzone)$, with 
impedances $\sqrt{\tau}$ and $-\sqrt{\tau}$, respectively.

\subsection{Finite volume discretization}
Having the complete eigenstructure of the hyperbolic part of the system,
it is straightforward to construct an explicit finite volume discretization,
using Strang splitting to handle the algebraic source term in a locally implicit manner.
\paragraph{Hyperbolic part.}
For the homogeneous part of the system,
i.e. the left-hand side of  \eqref{eq:ND-NLSH-real-imaginary},
we first consider the  Riemann problem in the $j$-th direction, given by
\begin{align}
    \label{eq:ND-NLSH-RP}
\partial_t \eta + B_j \partial_{x_j} \eta = 0, \qquad
\eta(\bfx,0) = \begin{cases}
\eta_{\ell} & \text{if } x_j < x_j^0, \\
\eta_{r} & \text{if } x_j > x_j^0.
\end{cases} 
\end{align}
The exact solution of the linear Riemann problem \eqref{eq:ND-NLSH-RP} is
\begin{align}
    \eta(\bfx,t) = \begin{cases}
        \eta_{\ell} & x_j < x_j^0 + t\lambda_1, \\
        \eta^* & x_j^0 + t\lambda_1 < x_j < x_j^0 + t\lambda_4, \\
        \eta_{r} & x_j > x_j^0 + t\lambda_4,
    \end{cases}
\end{align}
where $\eta^*$ satisfies the jump conditions
\begin{align}
    \label{eq:exact RS jth}
    \begin{cases}
        B_j (\eta^* - \eta_{\ell}) = \lambda_1 (\eta^* - \eta_{\ell}), \\
        B_j (\eta_{r}-\eta^*) = \lambda_4 (\eta^* - \eta_{r}),
    \end{cases}
    \Rightarrow
    \begin{cases}
        \eta^* - \eta_{\ell} = \alpha_1 v_1+ \alpha_2 v_2, \\
        \eta_{r}-\eta^* = \alpha_3 v_3 + \alpha_4 v_4.
    \end{cases}
    \Rightarrow
        \eta_{r} - \eta_{\ell} = V \boldsymbol{\alpha} \;,
\end{align}
and the coefficients $\alpha_i$ can be obtained by solving the linear system
$\tilde{V} \boldsymbol{\alpha} = P_j^T (\eta_{r} - \eta_{\ell})$,
since $P_j^T P_j = I_{4}$.

We use this solution to update the cell averages 
with LeVeque's second order wave-propagation
formulation of Godunov's method \cite{leveque2002finite}. 
In two spatial dimensions it reads
\begin{subequations} \label{eq:FV scheme}
\begin{align}
\eta_{ij}^{n+1}=\eta_{ij}&-\frac{\dt}{\dx}\left(B_1^{+}\Delta \eta_{i-1/2,j}+B_1^{-}\Delta \eta_{i+1/2,j}\right)
-\frac{\dt}{\dy}\left(B_2^{+}\Delta \eta_{i,j-1/2}+B_2^{-}\Delta \eta_{i,j+1/2}\right)\\
&-\frac{\dt}{\dx}\left(\tilde F_{i+1/2,j}-\tilde F_{i-1/2,j}\right)-\frac{\dt}{\dy}\left(\tilde G_{i,j+1/2}-\tilde G_{i,j-1/2}\right),
\end{align}
\end{subequations}
where the fluctuations $B_1^{\pm}\Delta \eta_{i\mp 1/2,j}$ and $B_2^{\pm}\Delta \eta_{i,j\mp 1/2}$
account for the net contribution of the Riemann problem solution in the 
$x$ and $y$ directions at the cell interfaces.
The terms $\tilde F_{i\pm 1/2,j}$ and $\tilde G_{i,j\pm 1/2}$ correspond to 
Lax-Wendroff-type second order corrections and transverse propagation of information 
at cell corners, which
allow us to increase the Courant number $\Courant$ close to unity in two dimensions.
In Section \ref{sec:tests} below, no limiting is applied
for the second order corrections and $\Courant=0.99$ is used for all
the numerical tests.
For more details on this standard discretization, and for the extension
to logically rectangular mapped grids (see Appendix \ref{sec:rotational-invariance}),
the reader is referred to \cite[Ch. 21, 23]{leveque2002finite} and \cite{calhoun2008logically}.

The finite volume scheme is implemented in Clawpack \cite{mandli2016clawpack},
which allows us to seamlessly use Berger-Olliger-Collela patch-based adaptive mesh refinement
with the wave-propagation formulation
\eqref{eq:FV scheme} \cite{berger1998amrwp}.
We adopt a refinement criterion based on differences of the modulus of the wavefunction $\rho = |q_0|$.
Specifically, a cell is flagged for refinement if the first-order
undivided difference with the nearest neighborhood exceeds a threshold:
\begin{align}
    \label{eq:refinement criterion}
    |\rho_{i,j}-\rho_{i\pm 1,j\pm 1}|>\refcrit.
\end{align}
For all the tests where AMR is used, we set $\refcrit=0.01$.

\paragraph{Source term.}
For the algebraic source term, we have the ODE system
$\partial_t \eta(\bfx,t) = S \eta(\bfx,t)$,
where the operator
\begin{align}
    S = \begin{bmatrix}
        0 & \mathbf{0}_{1\times N} & -\kappa (z_0^2 + w_0^2)-V_0(\bfx) & \mathbf{0}_{1\times N} \\
        \mathbf{0}_{N\times 1} & \mathbf{0}_{N\times N} & \mathbf{0}_{N\times 1} & -\tau^{-1} I_N \\
        \kappa (z_0^2 + w_0^2)+V_0(\bfx) & \mathbf{0}_{1\times N} & 0 & \mathbf{0}_{1\times N} \\
        \mathbf{0}_{N\times 1} & \tau^{-1} I_N & \mathbf{0}_{N\times 1} & \mathbf{0}_{N\times N}
    \end{bmatrix} \;,
\end{align}
is independent of $\eta$ since $\partial_t(z_0^2+w_0^2)$ is constant for this subproblem.
Setting $a(\bfx) = \left(\kappa  (z_0^2 + w_0^2)+V_0(\bfx)\right)$,
this ODE can be solved exactly by noting that
\begin{align}
    \exp(\Delta t S) = \begin{bmatrix}
        \cos\left(\Delta t a\right) & \mathbf{0}_{1\times N} & -\sin\left(\Delta t a\right) & \mathbf{0}_{1\times N} \\
        \mathbf{0}_{N\times 1} & \cos(\tau^{-1} \Delta t) I_N & \mathbf{0}_{N\times 1} & -\sin(\tau^{-1} \Delta t) I_N \\
        \sin\left(\Delta t a\right) & \mathbf{0}_{1\times N} & \cos\left(\Delta t a\right) & \mathbf{0}_{1\times N} \\
        \mathbf{0}_{N\times 1} & \sin(\tau^{-1} \Delta t) I_N & \mathbf{0}_{N\times 1} & \cos(\tau^{-1} \Delta t) I_N
    \end{bmatrix}.
\end{align}
The exact solution of this subproblem allows us to
handle the $\bigo(\tau^{-1})$  stiffness of the source term
while choosing the timestep only through the CFL condition
corresponding to the hyperbolic subproblem $\dt\propto \tau^{-1/2} \min\{\dx,\dy\}$.
Of course, this implies that explicit time stepping is advantageous rougly when
$ \tau^{1/2} \gtrapprox \max\{\dx,\dy\}$.
On the other hand, if $\tau^{1/2} \ll \min\{\dx,\dy\}$, then the time step restriction
becomes as severe as for an explicit scheme for the original NLS equation.

In our numerical examples we use a moderate value of $\tau$ and apply the splitting
scheme just described.
But if $\tau$ is very small then the stiffness of the source term leads
to an increase in the splitting error;  in this case an unsplit scheme might be preferable.
Some results for a fully explicit unsplit scheme are presented in Appendix \ref{sec:fully-explicit}.

\paragraph{Boundary conditions.} 
The boundary conditions for NLSH are imposed through ghost cells,
which are filled with the appropriate values depending on the type of boundary condition.
For Dirichlet boundary conditions, we enforce $q_0=0$ on the boundary.
For Neumann boundary conditions, it suffices to take $\bfqone\cdot \bfn = 0$ on the boundary since,
by the relaxation equation, this implies that $\partial_{\bfn} q_0 = 0$.
Given a state $\eta^{\text{in}}$ in the interior of the domain,
we can prescribe the state $\eta^{\text{g}}$ in the ghost cell such that the 
middle state
of the Riemann problem solution normal to the boundary satisfies the desired boundary condition.
By the rotational invariance of the system (see Appendix \ref{sec:rotational-invariance}), we have
\begin{align}
    z_0^* = \frac{1}{2}\left[(z_0^{\text{in}}+z_0^{\text{g}}) - \sqrt{\tau} (\bfwone^{\text{g}}\cdot \bfn-\bfwone^{\text{in}}\cdot \bfn)\right], \quad 
    \bfwone^*\cdot \bfn = \frac{1}{2}\left[(\bfwone^{\text{in}}\cdot \bfn+\bfwone^{\text{g}}\cdot \bfn) - \tau^{-1/2} (z_0^{\text{g}}-z_0^{\text{in}})\right],\\
    w_0^* = \frac{1}{2}\left[(w_0^{\text{in}}+w_0^{\text{g}}) + \sqrt{\tau} (\bfzone^{\text{g}}\cdot \bfn-\bfzone^{\text{in}}\cdot \bfn)\right], \quad
    \bfzone^*\cdot \bfn = \frac{1}{2}\left[(\bfzone^{\text{in}}\cdot \bfn+\bfzone^{\text{g}}\cdot \bfn) + \tau^{-1/2} (w_0^{\text{g}}-w_0^{\text{in}})\right].
\end{align}
Hence, for Dirichlet boundary conditions ($q_0^*=0$), it suffices to set $q_0^{\text{g}} = -q_0^{\text{in}}$ and $\bfqone^{\text{g}} = \bfqone^{\text{in}}$,
while for Neumann boundary conditions ($\bfqone^*\cdot \bfn = 0$), we set $q_0^{\text{g}} = q_0^{\text{in}}$ and $\bfqone^{\text{g}}\cdot \bfn = -\bfqone^{\text{in}}\cdot \bfn$.

\subsection{ImEx pseudopectral discretization}
In addition to the explicit finite volume method described in the previous section,
we also implemented a higher-dimensional version of the discretization method
introduced in \cite{biswas2025hyperbolic}.
This method is based on a pseudo-spectral discretization in space and a high-order, asymptotic-preserving ImEx RK scheme in time,
which allows us to study problems with smooth solutions, periodic boundary conditions, and small values of $\tau$.
With this method, mass is conserved to machine precision by using relaxation
for the dynamically chosen time step \cite{ketcheson2019relaxation}.
For the tests in which this method is used, we use ARS(4,4,3) \cite{ascher1997implicit} as the time integrator.

\section{Numerical experiments} \label{sec:tests}
In this section we present numerical experiments, using the methods described
in Section \ref{sec:FVM}.
In Section \ref{sec:verification} we present Fourier pseudospectral simulations that confirm the
theoretical results from this work.
In Section \ref{sec:FV-simulations}, we present finite volume simulations on mapped, adaptively
refined grids, that illustrate the usefulness of the NLSH system to study problems with sharp gradients and complex geometries.

\subsection{Verification of theoretical results} \label{sec:verification}
\subsubsection{Modulational instability of plane waves}
\label{sec:MI numerical}
We begin by verifying the stability of plane wave solutions of NLSH described in section \ref{sec:MI theoretical},
restricting ourselves to the one-dimensional case for simplicity.
First, we illustrate the effect of hyperbolization on the
shrinking of the unstable band for low-wavenumber perturbations.
For this we take $\kappa=1$, $p=3$, and $V_0(x)=0$,
and allow the relaxation parameter to take the values $\tau\in \{0.3, 0.4, 0.49, 0.51\}$.
We initialize the solution as a constant background with a small perturbation
given by \eqref{eq:sideband ansatz NLSH}
$q_0(x,0) = \sigma + \varepsilon(x,0)$,
$q_1(x,0) = \mathbf{\delta} (x,0)$,
where
$\sigma=1$,  $\varepsilon_+=10^{-3}$ and $k=1$.
The perturbation is chosen in the direction of the eigenvector corresponding to the root of the quartic \eqref{eq:ND-NLSH-k0-perturbation-quartic}
with the largest imaginary part.

\begin{figure}[]
    \centering
    \begin{subfigure}{0.49\textwidth}
        \includegraphics[width=\linewidth]{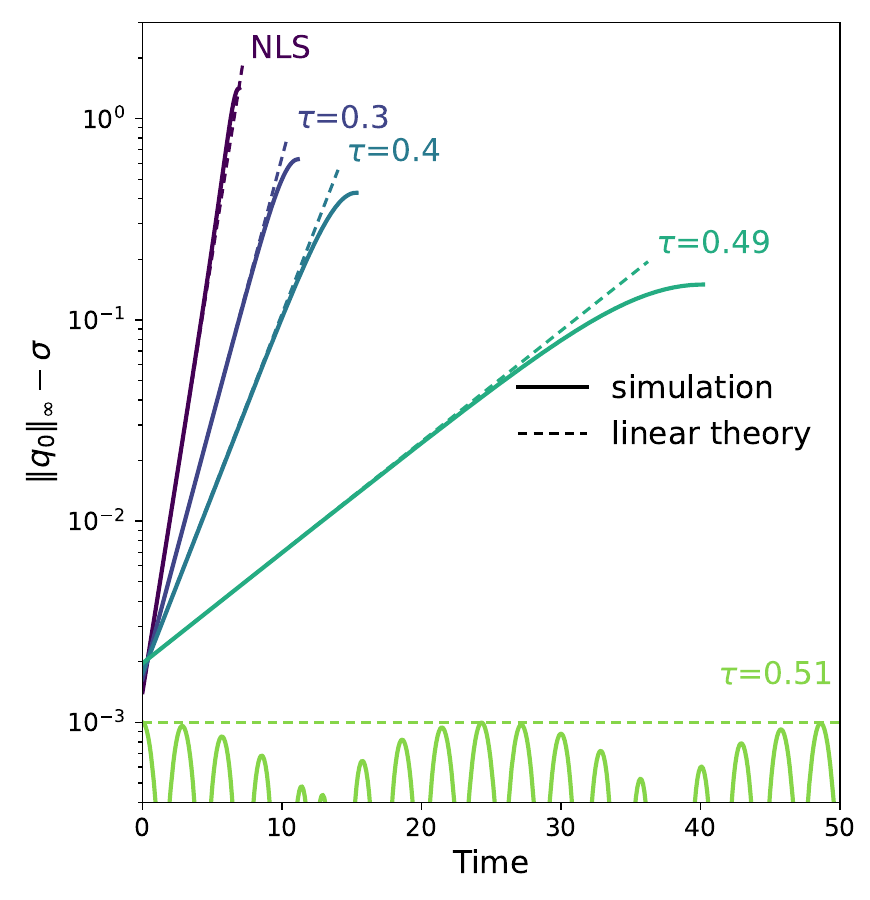}
    \end{subfigure}
    \begin{subfigure}{0.49\textwidth}
        \includegraphics[width=\linewidth]{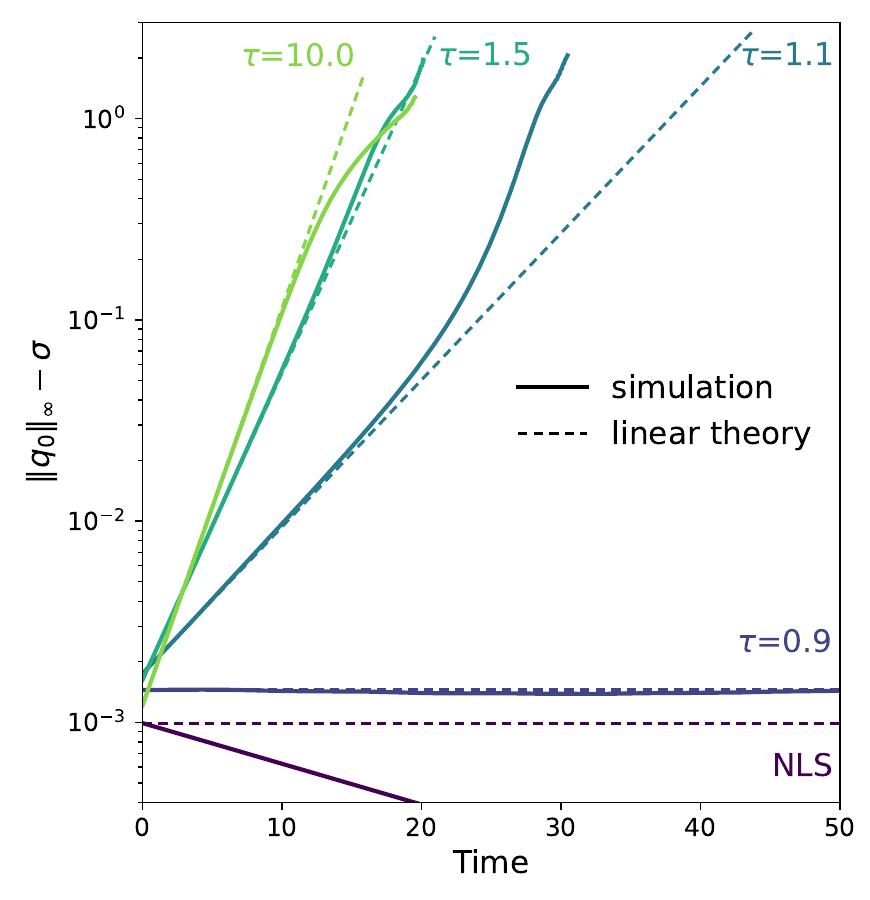}
    \end{subfigure}
    \caption{The deviation of $\|q_0\|_\infty$ from the background amplitude $\sigma$ 
    is plotted for different values of $\tau$.
        Left: Focusing low-wavenumber instability with is suppressed
        for $\tau>\tau_\text{low}= 1/2$ with $k=1$.
    Right: Defocusing high-wavenumber instability with is triggered
        for $\tau>\tau_\text{high}\approx 1$ with $k=8$.
    \label{fig: MI growth}}
\end{figure}

\begin{figure}[]
    \centering
    \begin{subfigure}{0.49\textwidth}
        \includegraphics[width=\linewidth]{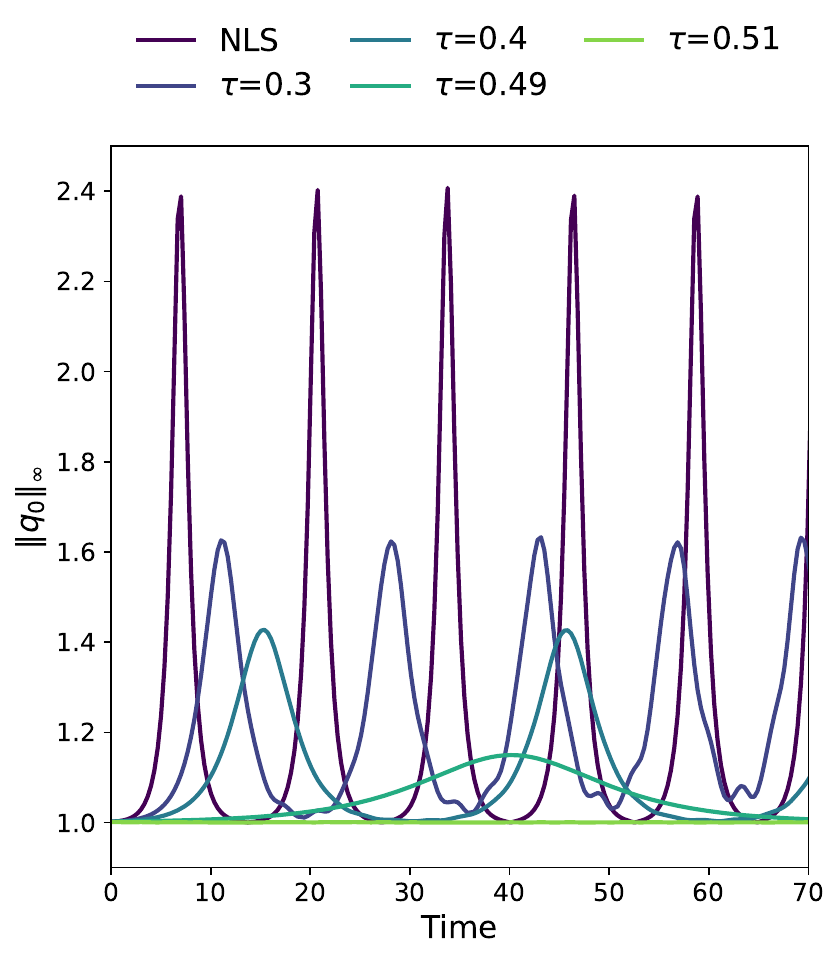}
    \end{subfigure}
    \begin{subfigure}{0.49\textwidth}
        \includegraphics[width=\linewidth]{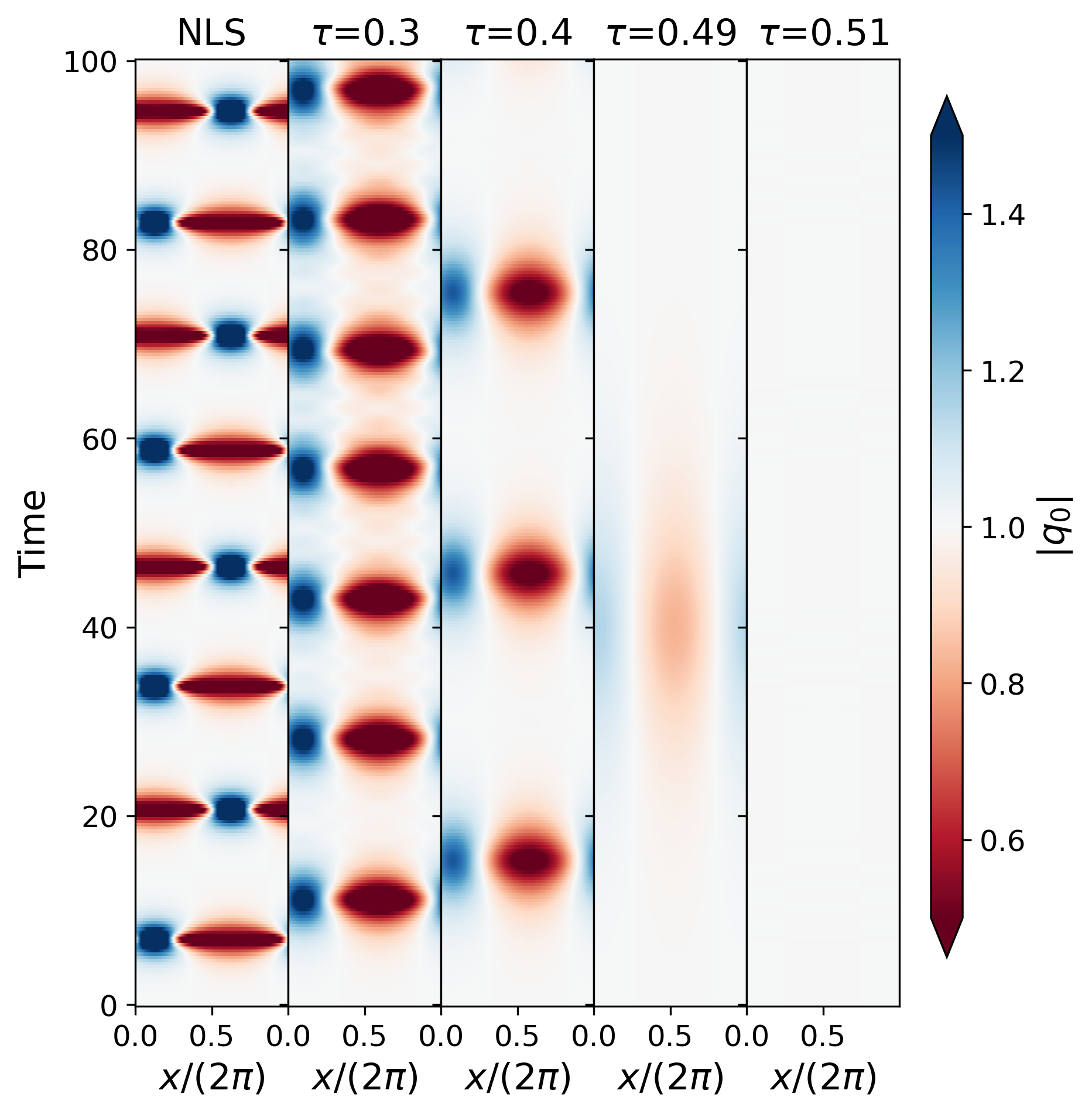}
        \label{fig: MI focusing final solution}
    \end{subfigure}
    \caption{Akhmediev breather-like solutions of the focusing NLSH system.
        Left: Time evolution of the maximum amplitude $\|q_0\|_\infty$.
    Right: Space-time contour plot of $|q_0|$ for the focusing NLSH system with 
    different values of $\tau$.
    \label{fig: MI focusing}}
\end{figure}

A domain of length $L=2\pi$ is used with periodic boundary conditions and a uniform grid of $N=2^{10}$ points.
The solution is evolved with the ImEx pseudospectral scheme until $T=100$.
In the left panel of Figure \ref{fig: MI growth}, we plot the deviation of the maximum amplitude $\|q_0\|_\infty$
from the background amplitude $\sigma$
for different values of $\tau$,
versus the value attained by the theoretical deviation $ \varepsilon_+ \exp (\max\{\mathrm{Im}(\Omega(\tau))\} t) - \sigma,$
predicted by the linear stability analysis in section \ref{sec:MI theoretical},
where the maximum is taken over roots of the quartic \eqref{eq:ND-NLSH-k0-perturbation-quartic}
for $k=1$.
The numerical results match the theoretical predictions very well:
the instability takes place when $\tau$ is chosen such that \eqref{eq: MI condition NLSH}
holds, which for the parameters in this experiment corresponds to
$\tau<\tau_{\text{low}}\coloneqq 1/2$.
Furthermore, as shown in Figure \ref{fig: MI focusing}, the amplitude of the perturbed solution does not grow unboundedly
when this instability takes place.
Instead, it grows until it is comparable to the background, then saturates and decreases toward the initial profile
periodically in time.
Moreover, the amplitude and frequency of the recurrence seems to decrease as $\tau$ increases.
This Fermi-Pasta-Ulam-Tsingou-like recurrence phenomenon is a characteristic feature of the focusing NLS equation \cite{yuen1978relationship},
and it suggests the existence of Akhmediev breather solutions of the NLSH system, whose study is left for future work.

Now we turn our attention to the high-wavenumber unstable band induced by the hyperbolization.
For this, we consider the defocusing NLSH system with $\kappa=-1$, $p=3$, $V_0(x)=0$,
and the same grid and initial conditions as in the previous test. 
We fix the wavenumber of the perturbation to $k=8$, and allow the relaxation parameter to take the values $\tau\in \{0.9, 1.1, 1.5, 10\}$.
For this instability to occur, the two conditions in \eqref{eq: MI condition NLSH D<0} must be satisfied,
which for this choice of parameters corresponds to $\tau>\tau_\text{high}\approx 1$.
Again, there is a remarkable agreement between the theoretical predicitions and the numerical results,
as shown in the right panel of Figure \ref{fig: MI growth}.
Contrary to the low-wavenumber instability, after the perturbation leaves the linear regime,
the solution does not return to the initial profile, but instead follows a more complex, turbulent, dynamics.
This can be explained by the fact that the instability has the effect of broadening the Fourier spectrum of the solution.
When this occurs in the case of the low-wavenumber instability, the new Fourier modes leave the unstable band and thus
the solution stabilizes.
On the other hand, in the case of the high-wavenumber instability the new Fourier modes remain in the unstable band and thus
a positive feedback loop is created.
This particularly destructive behavior further demonstrates the importance of a carefully chosen
relaxation
parameter, as discussed in Section \ref{sec:MI theoretical}.

\subsubsection{Well-prepared initial data and slow manifold convergence}
\label{sec:WP-IC numerical}
\begin{figure}
    \centering
    \includegraphics[width=\linewidth]{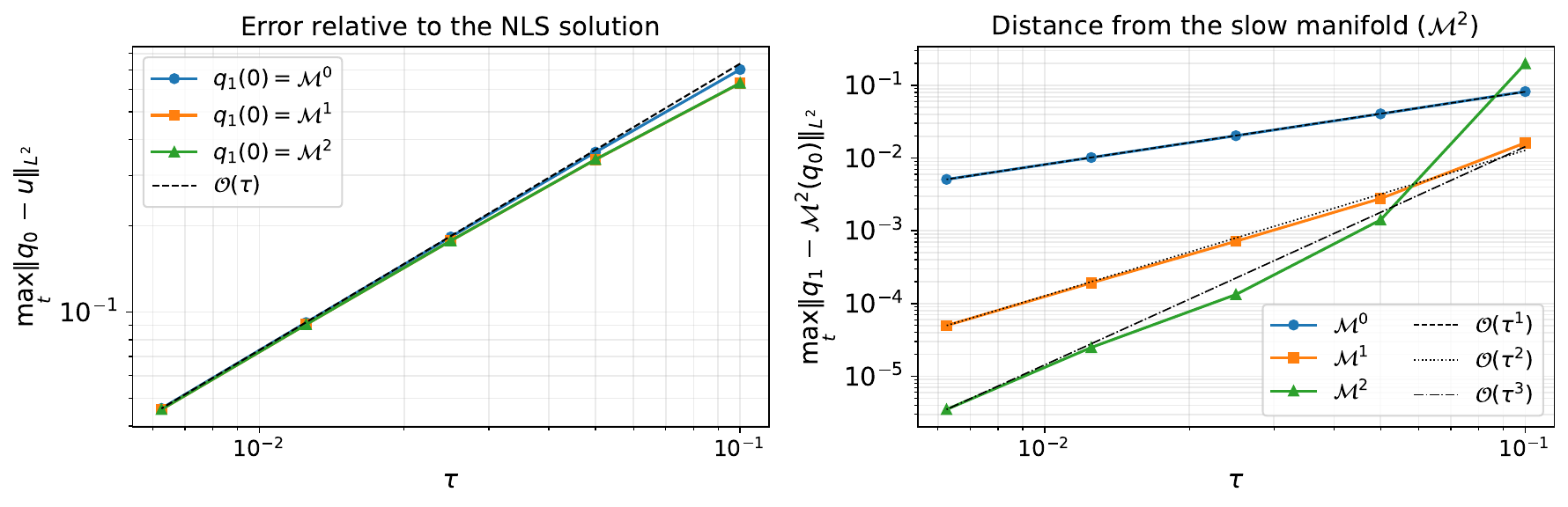}
    \caption{
        Convergence of NLSH solutions to NLS solutions as $\tau \to 0$ for different orders of well-prepared initial data.
        \label{fig: well-prepared convergence}}
\end{figure}

We now illustrate the effect of well-prepared initial data on the convergence of NLSH solutions to NLS solutions as $\tau \to 0$.
For this, we consider a one-dimensional ground state solution of NLS with $V_0(x)=0$, $\kappa=2$, $p=3$, and $\mu=1$, for which the exact solution is given by
$u(x) =  e^{it}\text{sech}(x)$.
For the values of $\tau \in \{10^{-1}, 5 \times 10^{-2},  2.5 \times 10^{-2}, 1.25 \times 10^{-2}, 6.25\times 10^{-3}\}$,
we initialize the NLSH system with well-prepared initial data $\bfqone(\cdot,0)=\smapiter{n}(\cdot,0)$ for $n=0,1,2$.
The solution is evolved until $T=5$ with the pseudo-spectral ImEx scheme 
($512$ grid points, $\Delta t =2\times 10^{-3}$),
and the error ($L^2$ in space, $L^\infty$ in time) between the NLSH solution and the exact NLS solution is computed.
Moreover, we also compute the difference between the auxiliary variable $q_1$ and the
highest available approximation of the slow manifold $\slowmap$, $\smapiter{2}$.
The results are shown in Figure \ref{fig: well-prepared convergence}.
As expected, the difference between the NLSH solution and the NLS solution decreases linearly with $\tau$ for all orders of well-prepared initial data,
as shown in the left panel of Figure \ref{fig: well-prepared convergence}.
Moreover, in the right panel of Figure \ref{fig: well-prepared convergence}, 
we can see that  higher-order well-prepared initial leads to the expected higher-order convergence of $\bfq$ to
the solution of the modified NLS equation, which is consistent with the formal analysis in Section \ref{sec: well-prepared initial data}.
These results suggest that the slow manifold of the NLSH system is $\bigo(\tau)$-close to the slow manifold of the NLS system, 
and that the rate of convergence for $q_0$ obtained in Theorem \ref{thm: convergence ND} is optimal.

\subsubsection{Convergence of ground state solutions}
\label{sec:numerical_ground_state}
\begin{figure}
    \centering
    \begin{subfigure}{0.49\textwidth}
       \includegraphics[width=\linewidth]{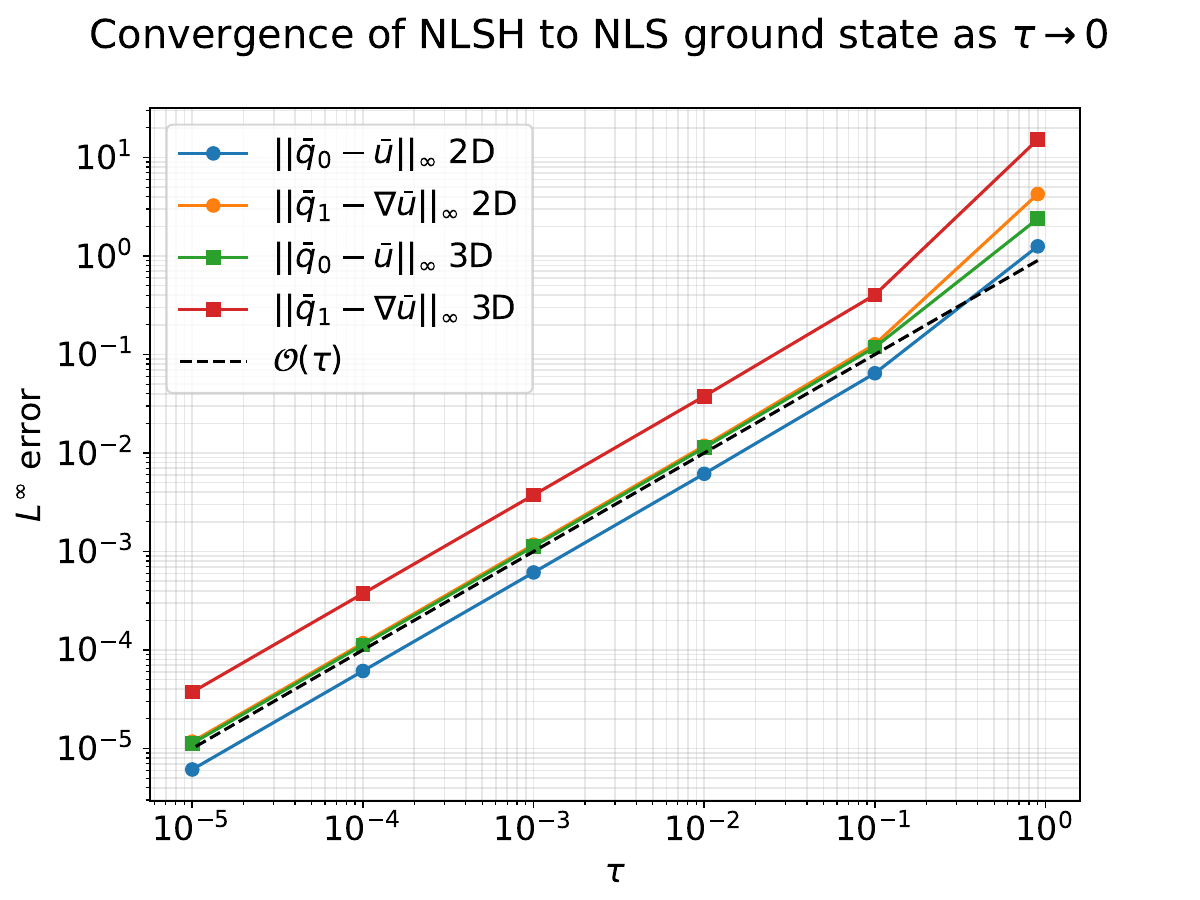}
	\caption{
        Convergence of NLSH to NLS ground states. Profiles generated with Petviashvili's method
        in 2D and 3D for different values of $\tau$.}
		\label{fig:ground_state_error_NLSH-NLS}
    \end{subfigure}
    \begin{subfigure}{0.49\textwidth}
        \includegraphics[width=\linewidth]{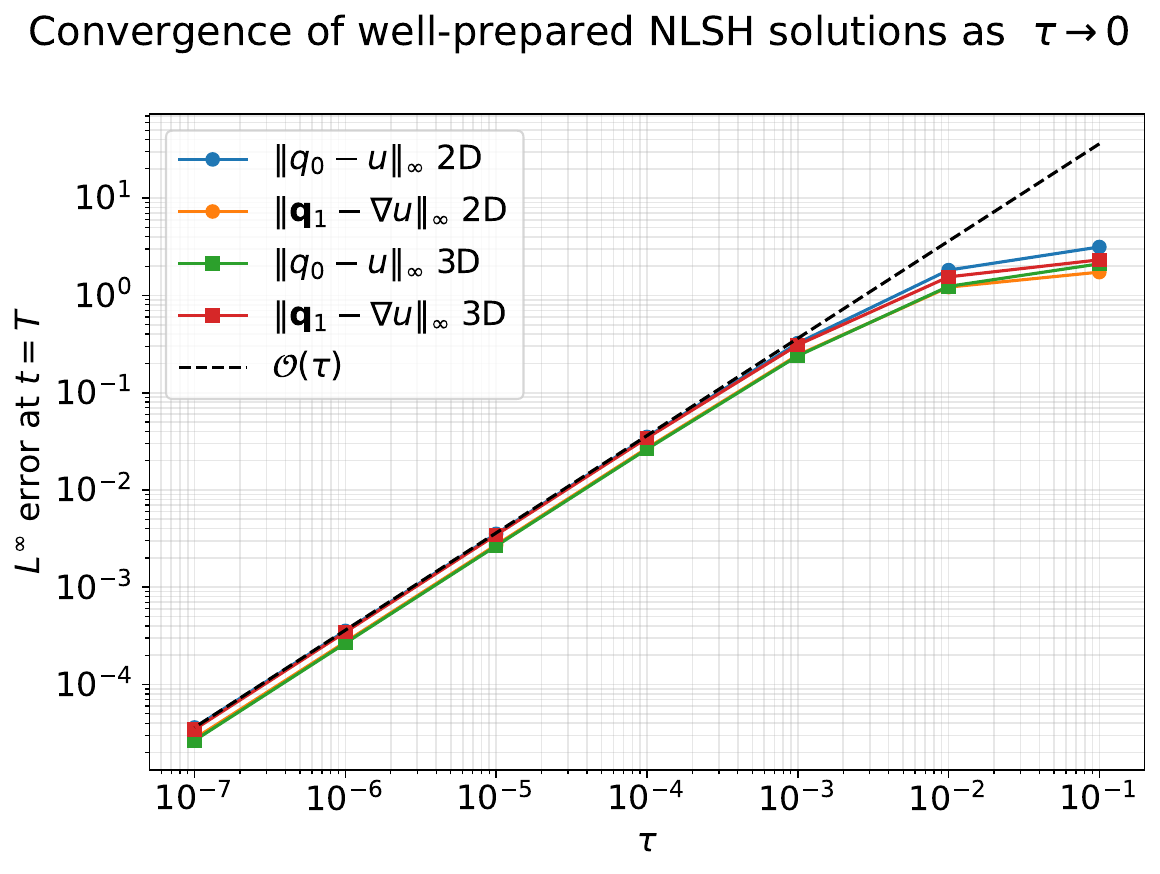}
	\caption{
        Convergence of NLSH to NLS solutions with well-prepared initial data based on NLS ground states.
        Initial data generated with Petviashvili's method.
    }
	\label{fig:ground_state_2d_err}
    \end{subfigure}
    \caption{Numerical verification of Proposition \ref{prop: converge ground states} and Theorem \ref{thm: convergence ND}
    using ground state solutions
    (see Section \ref{sec:numerical_ground_state})}
\end{figure}

Next, we validate the convergence results
presented in Proposition \ref{prop: converge ground states} and Theorem \ref{thm: convergence ND}.
For this, we consider NLS and NLSH with $V_0(\bfx)=0$, $\kappa=1$, and $p=3$ in two and three dimensions.
We generate a ground state solution of NLS using Petviashvili's method \cite{petviashvili1976equation}, \cite[Ch. 7]{yang2010nonlinear}
using a Fourier pseudo-spectral discretization with $512^2$ and $128^3$ grid points in 2D and 3D, respectively
with a periodic domain of length $L=20$ in each direction.
To validate Proposition \ref{prop: converge ground states}, we generate ground state solutions of NLSH for 
$\tau\in\{0.9, 10^{-1}, 10^{-2}, 10^{-3}, 10^{-4}, 10^{-5}\}$
using Petviashvili's method and compute the $L^\infty$ distance between $\qbar_0$, $\bar{\mathbf{q}}_1$ and $\bar{u}$, $\nabla \bar{u}$, respectively.
The results, shown in Figure \ref{fig:ground_state_error_NLSH-NLS}, confirm the
uniform, and linear in $\tau$, convergence of NLSH ground states to NLS ground states as $\tau \to 0$.

We also initialize the NLSH system with well-prepared (to order $0$, i.e. $\bfqone(\cdot,0)=\nabla q_0(\cdot,0)$) 
initial data based on the NLS ground state solution in two and three dimensions.
The solution is evolved with the pseudo-spectral ImEx scheme (ARS(4,4,3)) for
$\tau\in\{10^{-1},..., 10^{-7}\}$.
In the 2D case, we use $128^2$ grid points, a final time of $T=5$, and
a baseline time step of $\Delta t = 2\times 10^{-3}$,
which decreases like $\tau^{1/3}$ ($3$rd order time integration)
and vary slightly to conserve mass through relaxation.
In the 3D case, we use $48^3$ grid points, and a baseline time step of $\Delta t = 10^{-3}$.
Since the three-dimensional ground state is supercritical for $p=3$ and the solution is unstable,
we use a final time of $T=1$.
The results, shown in Figure \ref{fig:ground_state_2d_err}, confirm the convergence of NLSH solutions to NLS solutions as $\tau \to 0$
at rate $\bigo(\tau)$ for both $q_0$ and $\bfqone$. 
This is consistent with the results reported in \cite{giesselmann2025convergence}, and
suggests that the $\bigo(\sqrt{\tau})$ convergence estimate for $\bfqone$ might be suboptimal.

\subsection{Adaptively refined finite volume simulations} \label{sec:FV-simulations}

\subsubsection{Verification of the finite volume scheme}
\label{sec: verification FV}
\begin{figure}
    \centering
    \includegraphics[width=1\textwidth]{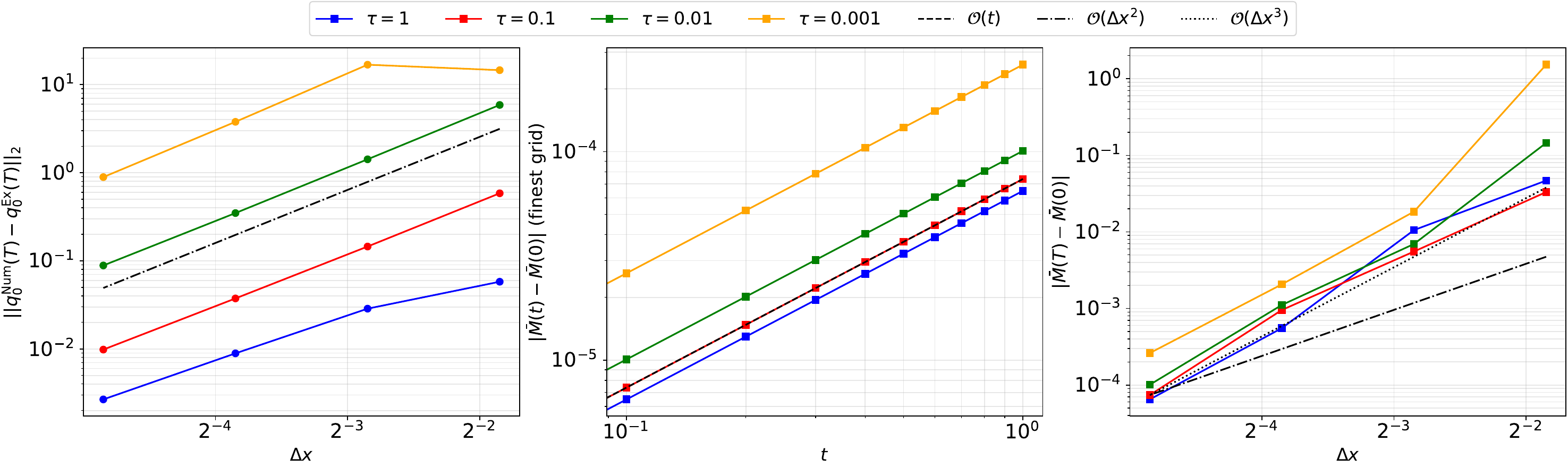}
    \caption{Convergence of the finite volume scheme for two-dimensional plane wave solutions of the NLSH system.
    The $L^2$ error is computed between the numerical solution and the exact solution for different grid resolutions
    (see Secction \ref{sec: verification FV} for details).}
    \label{fig: FV convergence OS}
\end{figure}

We now verify the convergence of the finite volume scheme described in Section \ref{sec:FVM} 
for two-dimensional plane wave solutions of the NLSH system \eqref{eq:ND-NLSH} with $V_0(\bfx)=0$, $\kappa=1$, and $p=3$.
The exact solution is given by the physical branch of \eqref{eq:ND-plane-wave-NLSH} with wavevector $\bfk=(1/\sqrt{2},1/\sqrt{2})$ and 
amplitude $\sigma=1$.
The computational domain is $[-\sqrt{2}\pi,\sqrt{2}\pi]^2$ with periodic boundary conditions, and the solution is evolved until $T=1$.
The convergence of the scheme is verified by computing the $L^2$ error between the numerical solution and the exact solution for different grid resolutions.
The results are shown in the leftmost panel of Figure \ref{fig: FV convergence OS}, where we can see that the scheme converges with second-order accuracy, as expected.
One can also observe that the constant of the leading error truncation term is proportional to $\tau^{-1}$,
which is a consequence of the fractional stepping.
In Appendix \ref{sec:fully-explicit}, we present results for a fully explicit scheme in which this issue is partially
mitigated, at the cost of a more severe time step restriction ($\dt\propto \tau^{-1} \dx$ rather than $\dt\propto \tau^{-1/2}  \dx$).
While the scheme is not designed to conserve the modified mass $\bar{M}$ exactly, as shown in the center
and right panels of Figure \ref{fig: FV convergence OS}, the mass is conserved to a very good approximation,
drifting linearly in time with a rate that decreases like $\bigo(\dx^3)$ as the grid is refined.

\subsubsection{Soliton reflection}
\label{sec:soliton_reflections}
We consider a test problem proposed by Katsaounis and Mitsotakis in \cite{katsaounis2018reflection}, 
in which a gray soliton is reflected from a vertical or slanted wall.
For this, we consider the defocusing NLSH system with $V_0(\bfx)=0$, $\kappa=-2$, and $p=3$.
The initial condition is a well-prepared planar gray soliton
\begin{align}
    \label{eq:grey soliton}
    u(x,y,t) =\eta[\cos \xi+i \sin \xi \tanh (\sin \xi \eta(-x+2 \eta \cos \xi t))] e^{-2 i \eta^2 t}
\end{align}
with $\eta=1$ and propagation speed $-2\eta\cos \xi=-\sqrt{2}$, i.e $\xi=\pi/4$.

Following \cite{katsaounis2018reflection}, we consider the domain $[-5,5]\times[-1,1]$ for the vertical wall test
and a trapezoidal domain with corners $(-8,-1)$, $(-7,1)$, $(8,1)$, and $(8,-1)$ for the slanted wall one.
In order to test the influence of grid distortion, we intentionally use a non-cartesian grid even for the 
vertical wall case.  The grid we use is
the disjoint union of five trapezoidal blocks, with the maximum angle between the vertical grid lines and the $y$-axis being $21.8^\circ$.
For the slanted wall test, we use a single trapezoidal block, obtained through a simple bilinear mapping of the unit square.
For the vertical wall test, we use a grid of $800\times 160$ cells and $\tau=10^{-2}$,
while for the slanted wall test we use a grid of $2560\times 320$ cells and $\tau=10^{-3}$.
Neumann boundary conditions are imposed on all boundaries.

\begin{figure}
    \centering
    \begin{subfigure}{0.49\textwidth}
        \includegraphics[width=\linewidth]{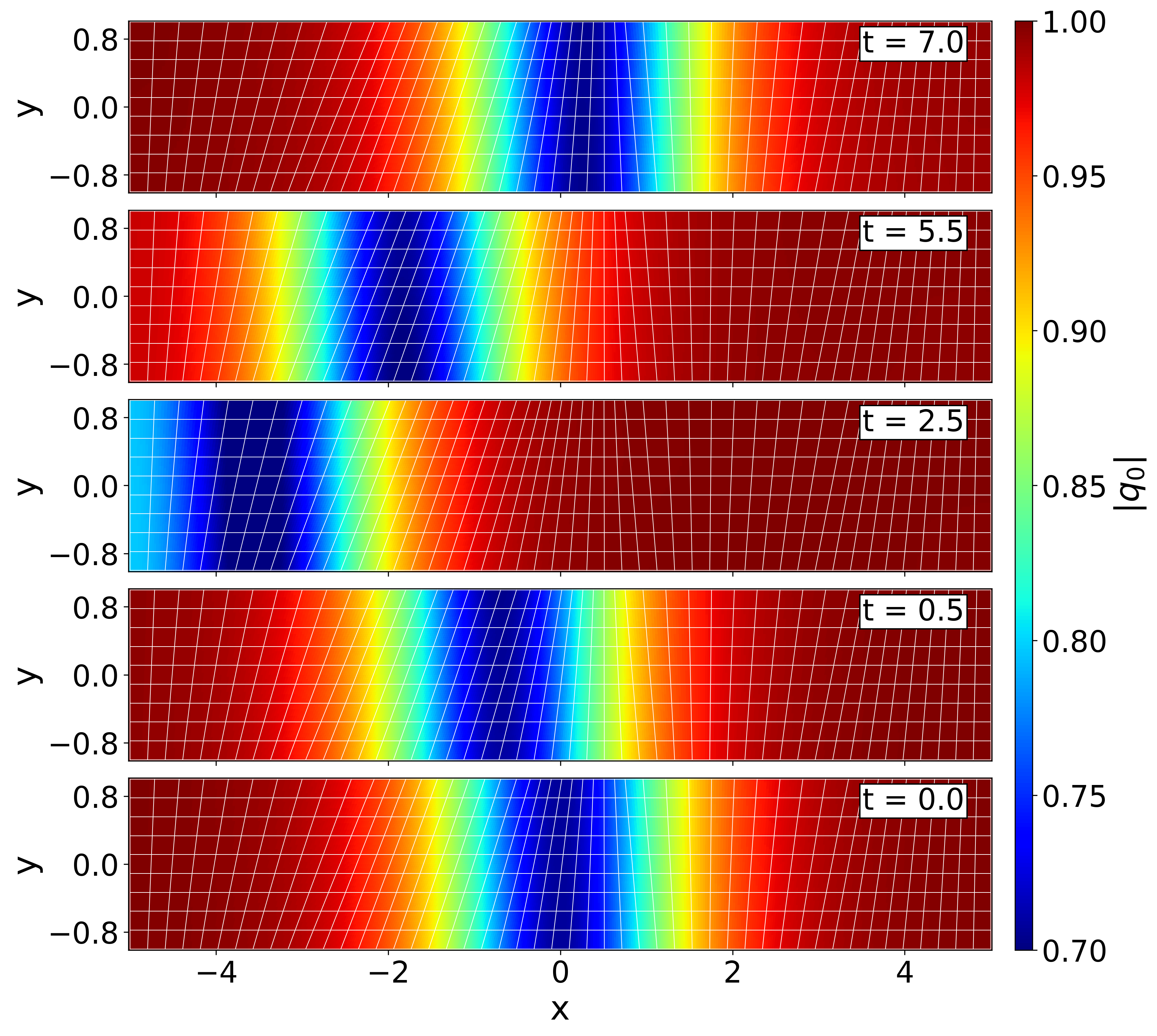}
        \label{fig:perfect_defocusing}
    \end{subfigure}
    \begin{subfigure}{0.49\textwidth}
        \includegraphics[width=\linewidth]{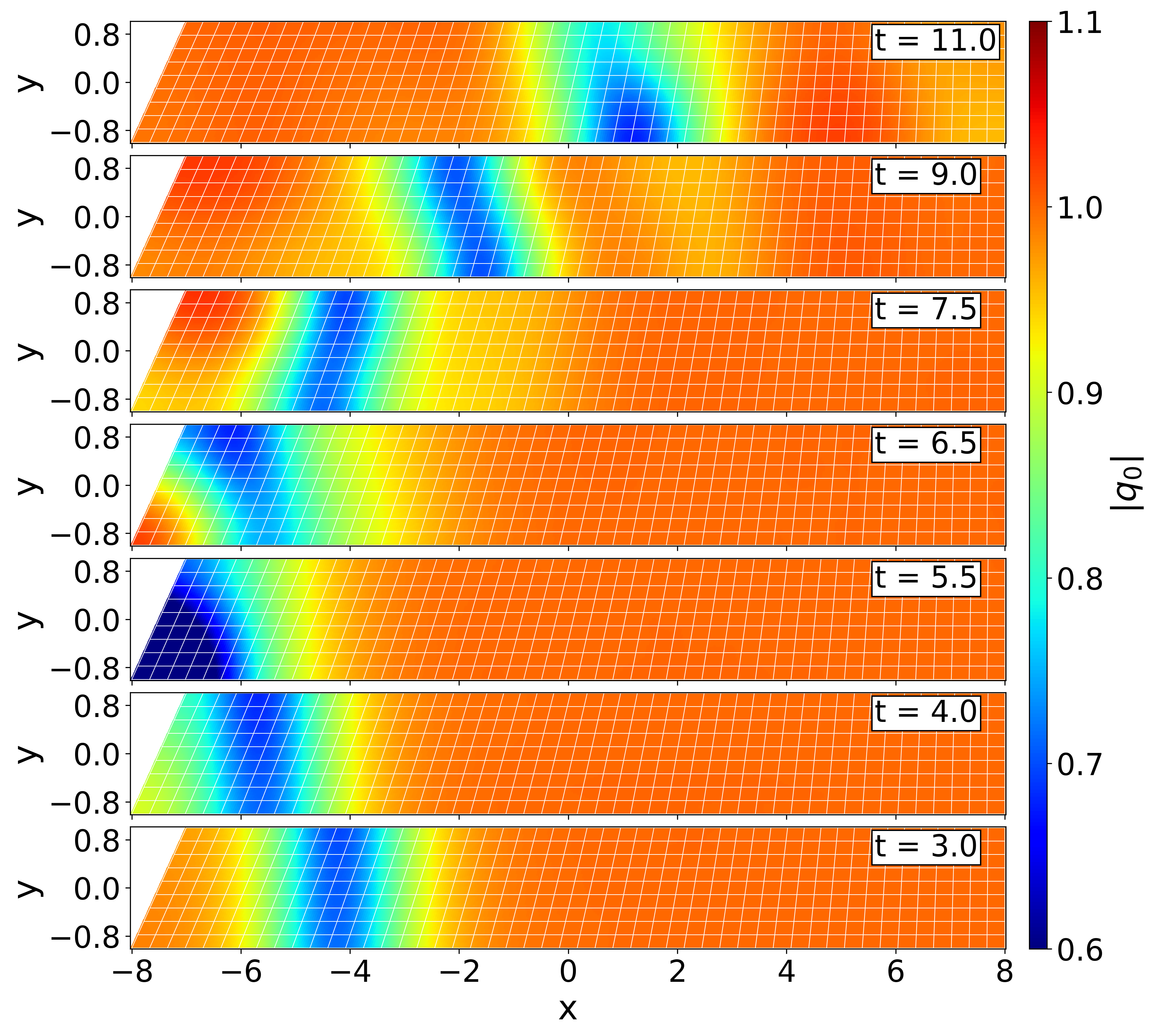}
        \label{fig:slanted_defocusing}
    \end{subfigure}
    \caption{NLS gray soliton reflections on solid walls with well-prepared initial data
    described in Section \ref{sec:soliton_reflections}.
    Undersampled grid lines displayed for reference.
    Left: perfect reflection on an orthogonal wall with an anisotropic grid.
    Right: reflection on a diagonal wall with a trapezoidal grid.
    \label{fig:soliton_reflections_defocusing}}
\end{figure}

Snapshots of the solution at different times are shown in Figure \ref{fig:soliton_reflections_defocusing},
along with coarse grid lines for reference.
As expected, the soliton is reflected off the vertical wall without any visible change in shape.
On the other hand, the structure of the soliton is distorted after reflecting off the slanted wall.
In both scenarios, excellent agreement is observed with respect to the reference results 
for NLS presented in \cite{katsaounis2018reflection} (see Figures 4 and 5 in that work).

\subsubsection{Transverse instability of bright solitons}
\label{sec:transverse instabilities}
It is well-known that plane-wave soliton solutions of the focusing NLS equation
are unstable to transverse perturbations \cite{yang2010nonlinear}. 
Here we study the transverse stability of soliton solutions of NLSH.
We consider a standing bright soliton solution of the one-dimensional focusing NLS equation
in the $x$-direction, given by 
\begin{align}
    u(x,t) = \sqrt{2} \, \mathrm{sech}(x) \, e^{i t}.
\end{align}
We set $q_0(x,y,0) = u(x,0)$, $\bfqone(x,y,0) = [\partial_x u(x,0), 0]^T$,
and add a small amount of random noise (1\% of the amplitude) to $q_0$ at every grid point
to trigger the transverse instability (following \cite[Ch. 5.9.3]{yang2010nonlinear}).
We consider the computational domain $[-20,20]\times[-20,20]$ with periodic boundary conditions in the $x$-direction
and reflecting boundary conditions in the $y$-direction.
Three levels of refinement are used, with the coarsest level covering the entire domain
with $80\times 80$ cells, and refinement factor of 5 in each direction between levels.

\begin{figure}
    \centering
    \includegraphics[width=0.9\linewidth]{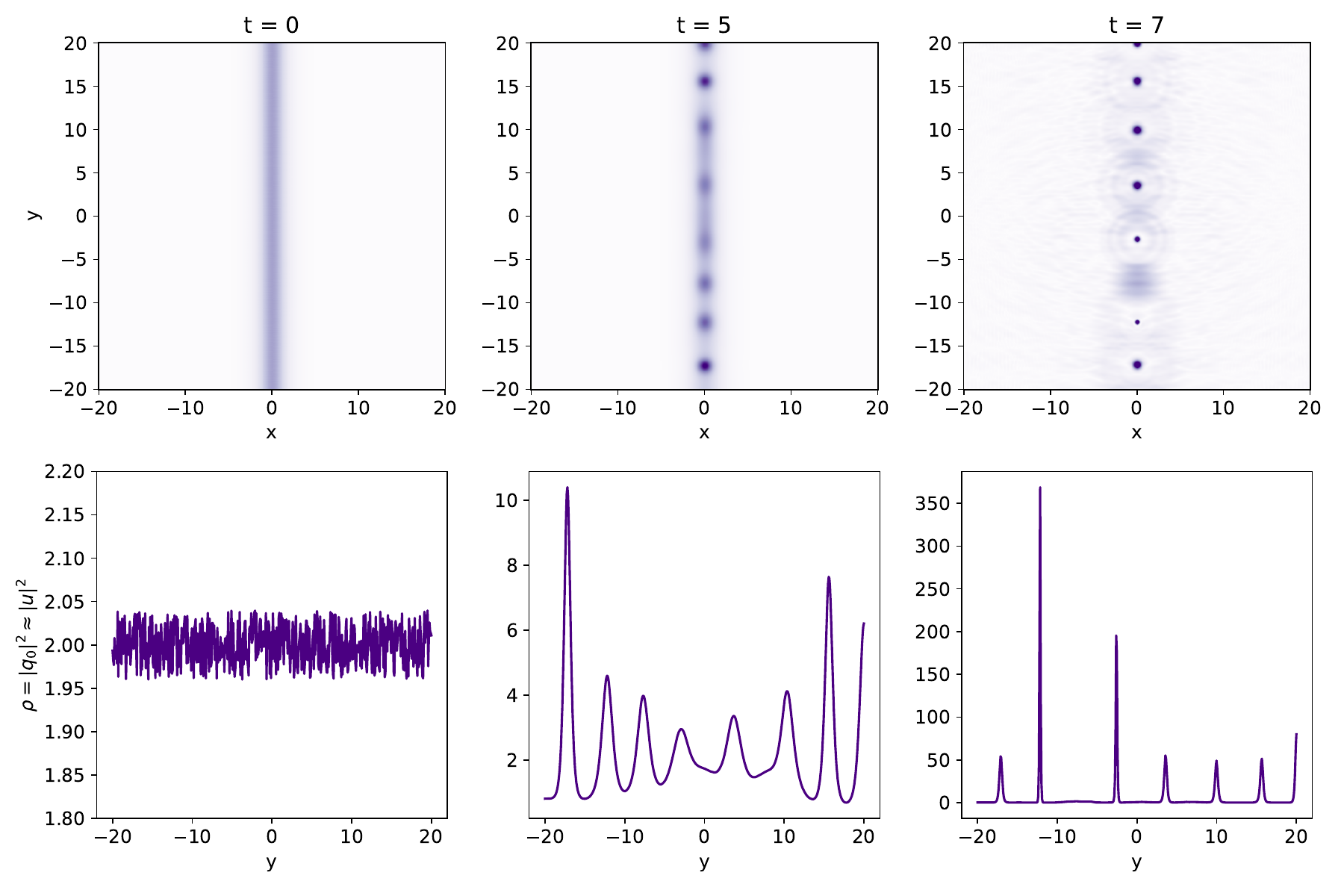}
    \caption{Transverse instability of a bright soliton in the focusing NLS equation.
    In the top row, we show the time evolution of the mass density $|q_0|^2$ at different times.
    The bottom row depicts a vertical slice of the mass density at $x=0$ at different times, where the development of neck instabilities can be observed.}
    \label{fig:transverse_instability_bright_soliton}
\end{figure}

As can be seen in Figure \ref{fig:transverse_instability_bright_soliton},
neck instabilities develop in the solution, and the soliton eventually breaks up into a series of filaments.
While the discretization used herein is able to capture
the early onset of this well-known phenomenon, it is important to note that it is not designed
to accurately resolve the blow-up of the solution for arbitrarily long times, 
which is a challenging problem in its own right that requires specialized numerical techniques \cite{akrivis2003numerical}.

\subsubsection{Radially symmetric dispersive shock waves}
\label{sec:radial DSW}
\begin{figure}[]
    \centering
    \begin{subfigure}{0.3\textwidth}
        \includegraphics[width=\linewidth]{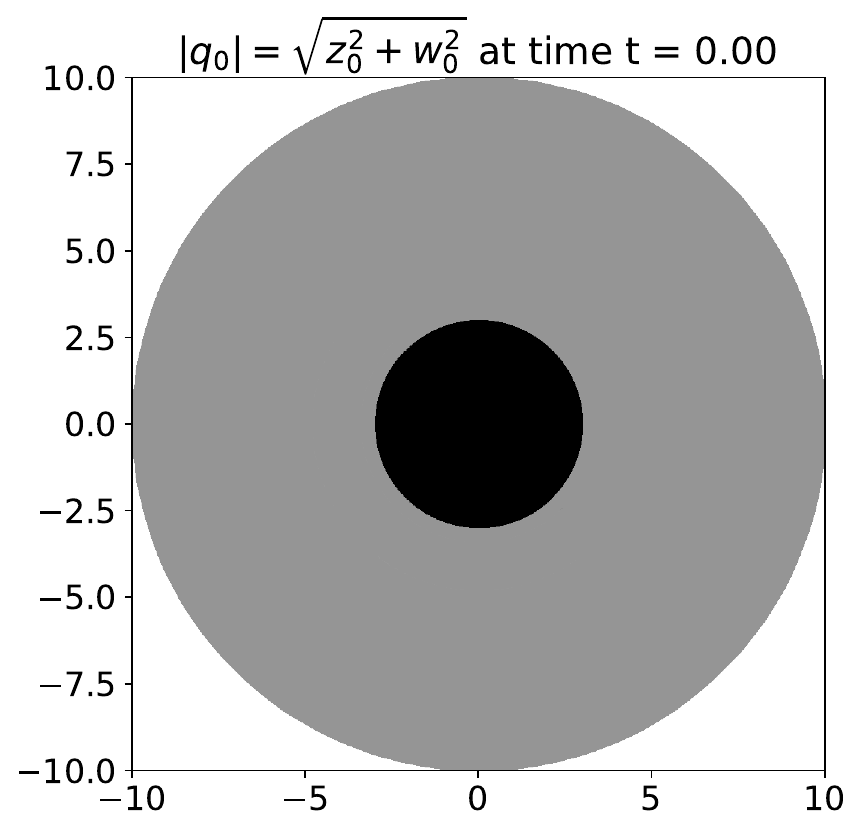}
        \caption{}
        \label{fig:dsw_t0}
    \end{subfigure}
    \begin{subfigure}{0.3\textwidth}
        \includegraphics[width=\linewidth]{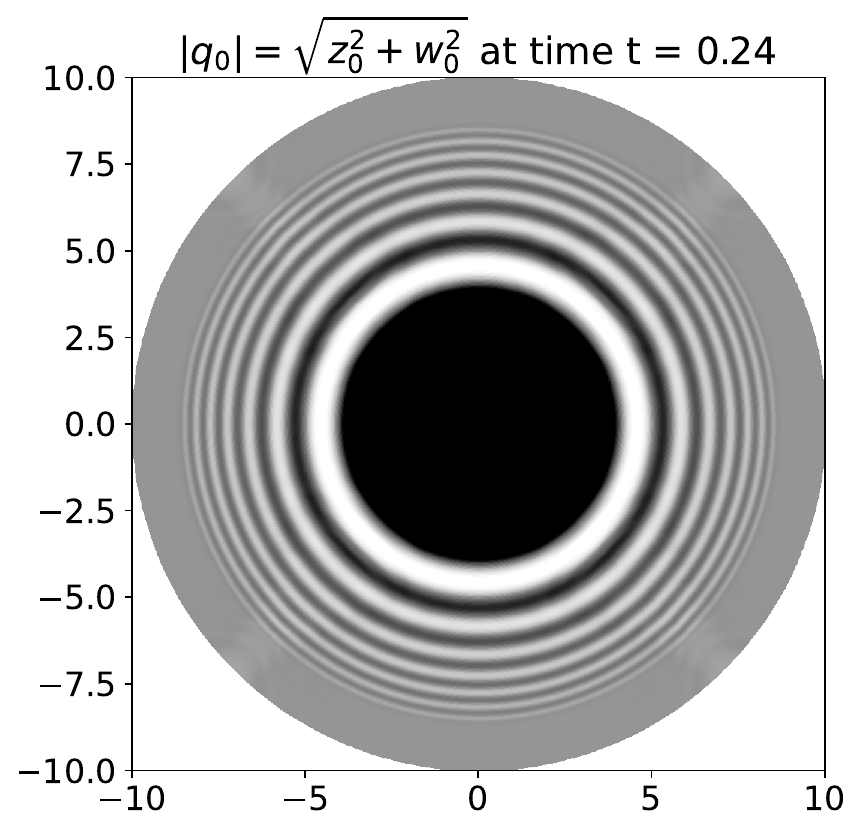}
        \caption{}
        \label{fig:dsw_t2}
    \end{subfigure}
    \begin{subfigure}{0.27\textwidth}
        \includegraphics[width=\linewidth]{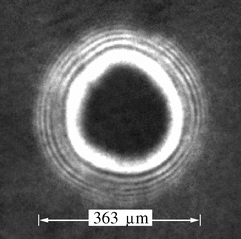}
        \caption{}
        \label{fig:Experiment_Hoefer}
    \end{subfigure}
    \vspace{1em} 

    \begin{subfigure}{0.3\textwidth}
        \includegraphics[width=\linewidth]{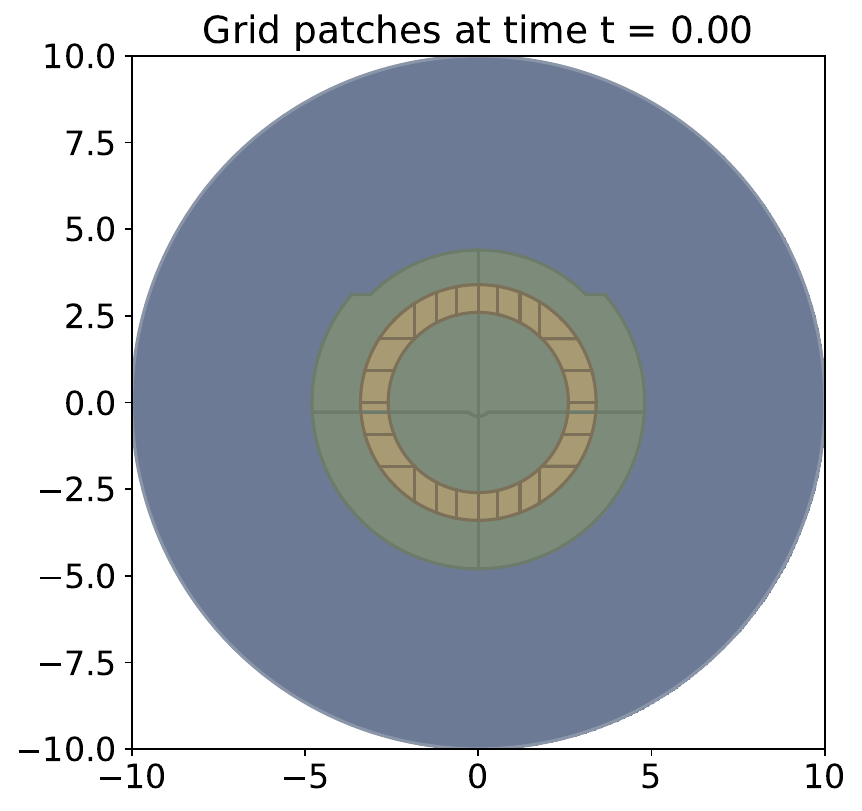}
        \caption{}
        \label{fig:dsw_t3}
    \end{subfigure}
    \begin{subfigure}{0.3\textwidth}
        \includegraphics[width=\linewidth]{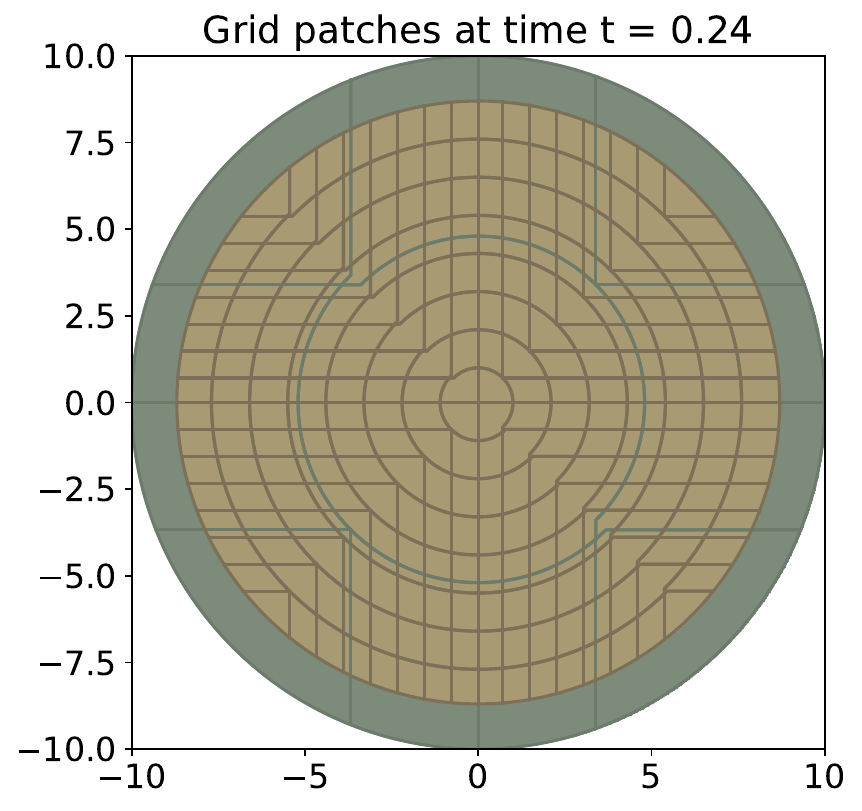}
        \caption{}
        \label{fig:dsw_t5}
    \end{subfigure}
    \begin{subfigure}{0.3\textwidth}
        \includegraphics[width=\linewidth]{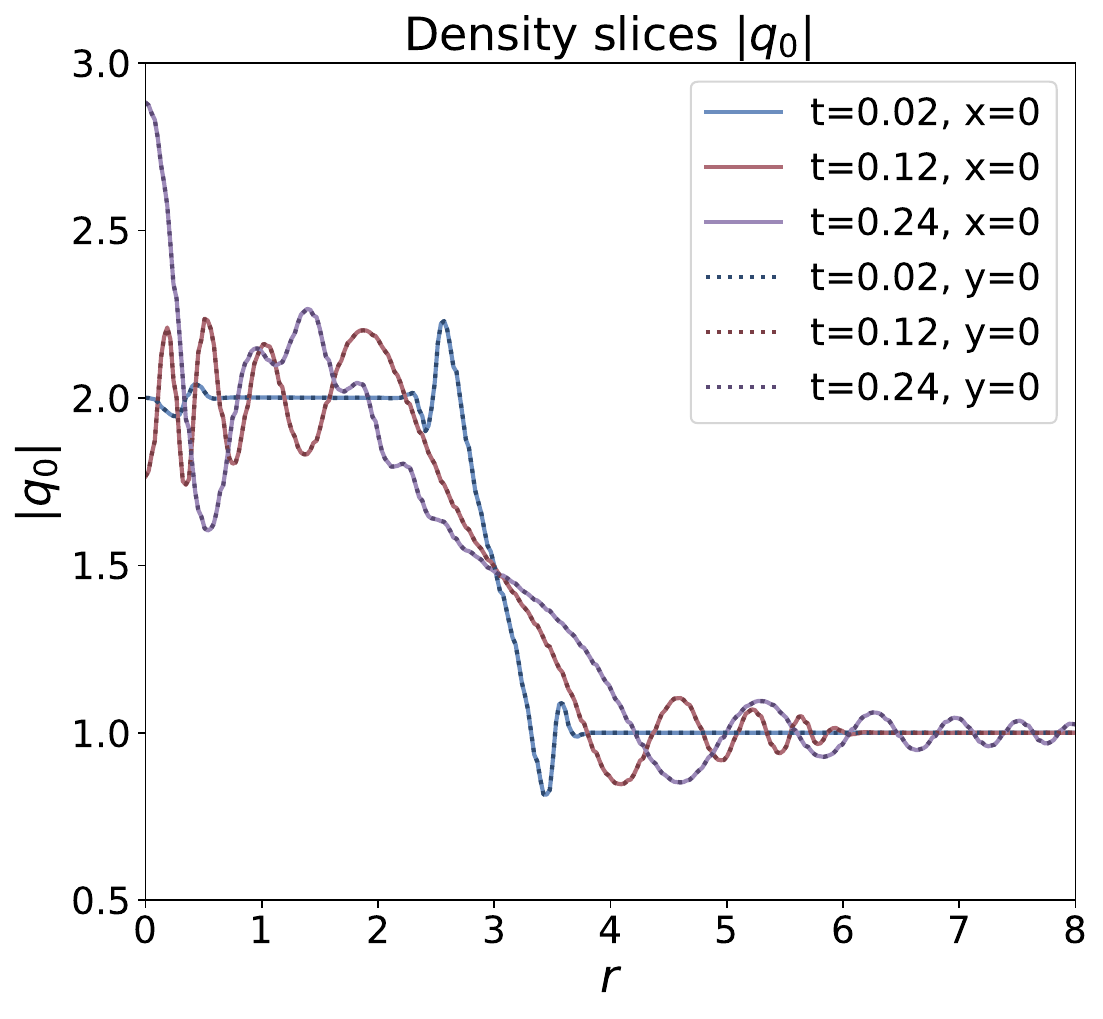}
        \caption{}
        \label{fig:dsw_slices}
    \end{subfigure}

    \caption{Radially symmetric dispersive shock wave on the disk with three levels of refinement.
    Top row: (a) and (b) depict the time evolution of mass density together with experimental data
     from \cite{hoefer2006dispersive} (c). Bottom row: refinement levels 1 (blue), 2 (green), 3 (golden), 
    and patch boundaries at initial (d) and final (e) times.
    (f) shows radial slices of the mass density at $y=0$ and $x=0$ at different times.
    Figure (c) is reproduced from \cite[Fig.~6]{hoefer2006dispersive}, with permission.}
    \label{fig:DSW_2D}
\end{figure}

We consider the defocusing NLSH equation ($\kappa=-1$)
in two spatial dimensions with a radially symmetric initial condition
given by the cylindrical Riemann problem:
\begin{align}
    q_0(r,0) = \begin{cases}
    \sqrt{\rho_i}, & r \leq r_0, \\
    \sqrt{\rho_o}, & r > r_0,
    \end{cases} \qquad
    \bfqone(r,0) = \mathbf{0},
\end{align}
where $r=\sqrt{x^2+y^2}$ is the radial coordinate, $r_0$ is the initial radius of the discontinuity,
and $\rho_i > \rho_o$ are the inner and outer densities, respectively.
This initial condition gives rise to a radially symmetric dispersive shock wave
that propagates outward from the origin. 
Proceeding similarly to \cite{busto2021high}, we set $\rho_i=2$, $\rho_o=1$, and $r_0=3$.
The computational domain is given by the square $[-1,1]\times[-1,1]$,
which is mapped to a disk of radius $R=10$ using the mapping proposed in \cite[Sec 3.2]{calhoun2008logically}.
Three levels of refinement are used, with the coarsest level covering the entire domain with $50\times 50$ cells,
and each finer level refining by a factor of 4 in each direction.
The relaxation parameter is set to $\tau=10^{-3}$, and
constant extrapolation boundary conditions are applied at the outer boundary of the disk.
Figure \ref{fig:DSW_2D} shows the time evolution of the mass density
and grids at initial and final times,
along with experimental data from \cite{hoefer2006dispersive} and radial slices of the mass density at different times.
Remarkable agreement is observed between the numerical results and the experimental data, 
which further validates the numerical results
obtained with a different system and discretization in \cite{busto2021high}.

\subsubsection{Flow past a cylinder}
\label{sec:flow past cylinder}
To leverage the advantages of mapped grids, we simulate the flow
of a superfluid
past an impenetrable
cylinder.
Following the setup in \cite{frisch1992transition}, we consider the defocusing NLSH system with $\kappa=-1$, $p=3$.
We use a plane wave of the form \eqref{eq:ND-plane-wave-NLSH} in the $x$ direction, with amplitude $\sigma=1$ and wavenumber $k=\frac{1}{2}$.
The relaxation parameter is set to $\tau=10^{-3}$.

\begin{figure}
    \centering
    \includegraphics[width=0.9\linewidth]{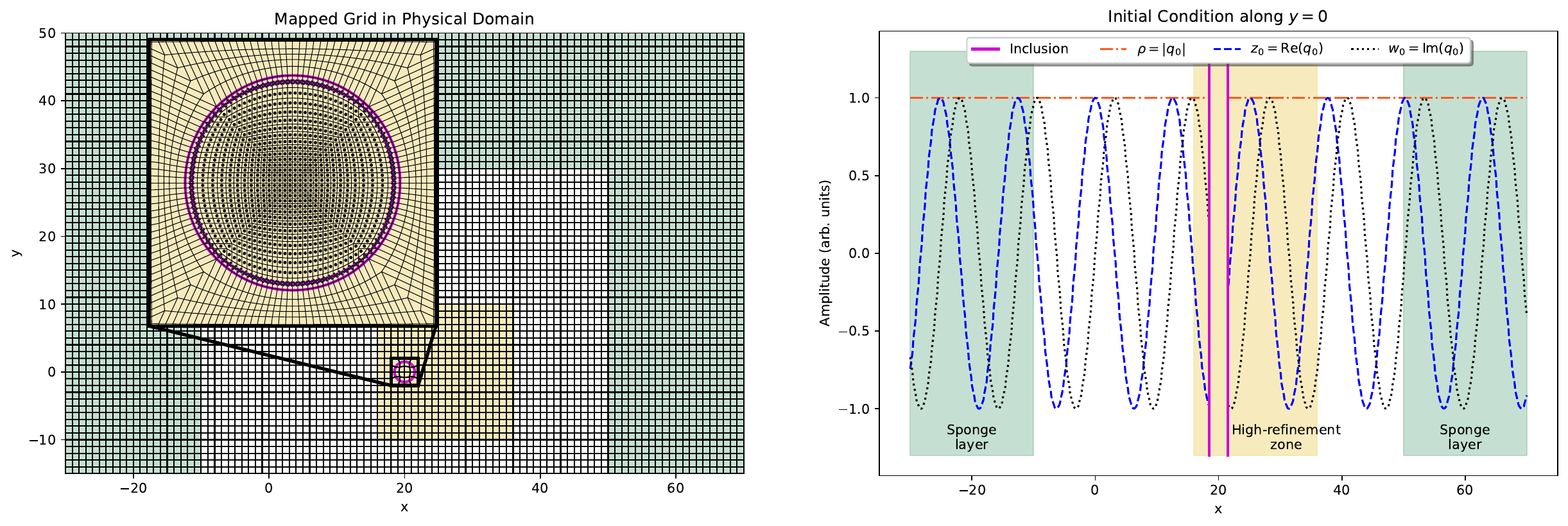}
    \caption{Setup for flow past a cylinder. Left: Mapped grid in physical domain. Right: Initial condition and sponge layers at the slice $y=0$.
    }
    \label{fig:flow_past_cylinder_setup}
\end{figure}

The computational setup is illustrated in Figure \ref{fig:flow_past_cylinder_setup} and described below.
The obstacle is centered at $(20,0)$ and has radius $r=3/2$, while the physical domain is given by the square $[-30,70]\times[-50,50]$.
Three refinement levels are used, with the coarsest level covering the entire domain with $200\times 200$ cells,
and each finer level refining by a factor of 5 in each direction.
Sponge layers of 20 units of width ($\approx 3/2$ of a characteristic wavelength) are used at the outer boundaries (green region)
to impose the plane wave inflow
and absorb the initial transient waves while minimizing reflections.
These layers are implemented simply by computing a weighted average of the solution with the desired plane wave solution at each time step \cite{larsen1983open},
and they are always forced to the second level of refinement.
The region $[16,36]\times[-10,10]$ (yellow region) around the cylinder is fixed at the finest level of refinement to better capture the flow features in the wake of the obstacle.
The remaining portion of the domain (white region)  is allowed to adapt dynamically based on the criterion \eqref{eq:refinement criterion}, with a maximum refinement level of 2.

\begin{figure}
    \centering
    \includegraphics[width=\linewidth]{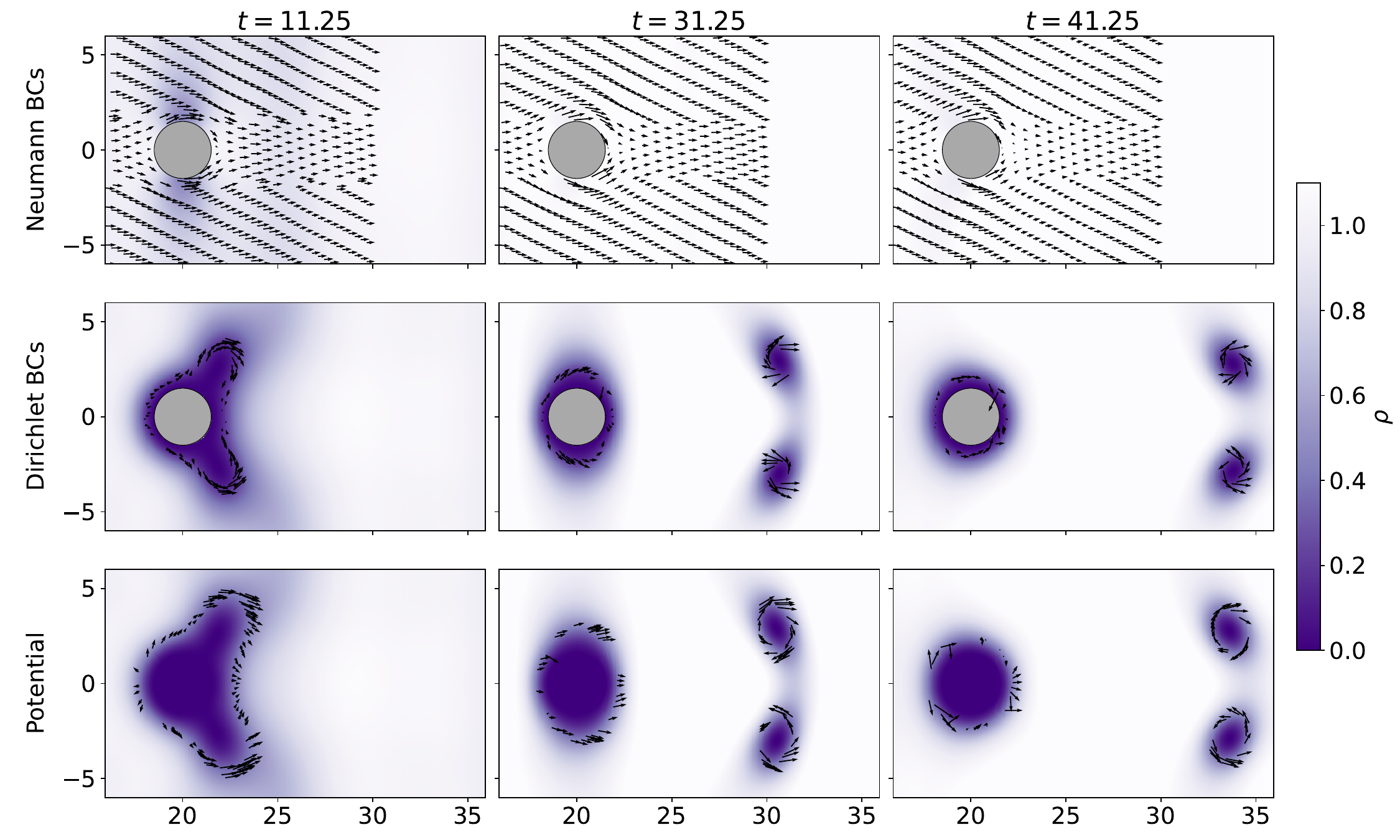}
    \caption{Superfluid flow past a cylindrical obstacle with different boundary conditions
    with a finite volume discretization of NLSH ($\tau=10^{-3}$) on a mapped grid (see Section \ref{sec:flow past cylinder})}
    \label{fig:flow_past_cylinder_results}
\end{figure}

We use a mapped grid \cite{calhoun2008logically} which allows us to conform to the shape of the cylinder, 
as illustrated in the left panel of Figure \ref{fig:flow_past_cylinder_setup}, thus
avoiding staircase effects.
The first-order formulation allows us to handle the boundary of the cylinder in a straightforward manner
(see Section \ref{sec:FVM}).
As in \cite{frisch1992transition}, we impose zero Dirichlet boundary conditions on the cylinder.
For the sake of comparison, we also performed simulations with homogeneous Neumann boundary conditions
and representing the cylinder through a large discontinuity in the potential:
\begin{align}
    V_0(\bfx) = \begin{cases}
        -10^{4}, & \text{if } |\bfx-(20,0)|<3/2, \\
        0, & \text{otherwise}.
    \end{cases}
\end{align}
Figure \ref{fig:flow_past_cylinder_results} shows the mass density and a quiver plot of the superfluid velocity,
given by
\begin{align}
    \nabla \theta(\bfx,t)=\frac{\mathrm{Im}\left(q_0^* \nabla q_0 \right)}{|q_0|^2} \approx \frac{\mathrm{Im}\left(q_0^* \bfqone\right)}{|q_0|^2}, 
\end{align}
where $\theta$ is the phase of the wavefunction $q_0 = |q_0| e^{i \theta}$.
We can observe the formation of vortices in the wake of the cylinder in the
case of Dirichlet boundary conditions and a repulsive potential, which is consistent with the results in \cite{frisch1992transition}.
As one would expect, the use of homogeneous Neumann boundary conditions suppresses the
formation of vortices, and the flow remains laminar.

\section{Conclusion and future work}
\label{sec: conclusion} 
We have proposed and studied NLSH \eqref{eq:ND-NLSH}, a first-order hyperbolic
relaxation of the nonlinear Schr\"odinger equation \eqref{eq:ND-NLS} in an
arbitrary number of spatial dimensions, extending the one-dimensional model of
\cite{ketcheson2025approximation,biswas2025hyperbolic}.  The system is
Hermitian, hyperbolic, and it
retains much of the structure of the original model.  
In particular, it is a
Hamiltonian system and it exactly
conserves modified mass, momentum, and energy functionals that reduce to the
corresponding NLS invariants in the relaxed limit.
 
We also examined several families of special solutions.  NLSH admits plane waves
whose physical dispersion branch  converges to that of NLS.
A modulational stability analysis of these plane waves shows
that hyperbolization acts on the instability structure in two distinct ways: it
shrinks the focusing low-wavenumber Benjamin--Feir band, recovering the NLS
condition in the relaxed limit, and it introduces a high-wavenumber unstable band
\eqref{eq: MI condition NLSH D<0} with no counterpart in NLS.  The latter is
absent provided the relaxation parameter is taken below an explicit threshold
that depends only on the nonlinearity strength.  
Standing waves behave similarly: the ground state of NLSH is a spatial rescaling
of the NLS ground state, and  converges to
it uniformly and linearly in the relaxation parameter.
 
We proved in Theorem \ref{thm: convergence ND} that
smooth enough solutions of NLSH  
converge to those of NLS over finite
time intervals. 
The proof  uses a relative energy argument in the spirit of
\cite{giesselmann2025convergence} and the modified mass as the stability
functional,
and
relies on using well-prepared initial data for the auxiliary variables.
The theorem further provides a linear convergence rate in the relaxation parameter, 
which is verified numerically for ground state solutions.
 
On the numerical side, we presented an explicit finite volume discretization of
NLSH that is compatible with mapped grids, embedded boundaries, and
patch-based adaptive mesh refinement, making it suitable for
scenarios where standard pseudospectral methods are awkward to apply. 
The solver is validated against gray soliton reflections, transverse instabilities
of bright solitons, radially symmetric dispersive shock waves, and vortex shedding past a cylinder
with moderate values of the relaxation parameter.  The results are consistent with the
expected NLS dynamics.
Moreover, for tests where the solution is smooth and periodic, we also implemented a 
pseudo-spectral discretization in space with a high-order, asymptotic-preserving ImEx RK scheme in time,
which allows for careful study of problems in which the relaxation parameter vanishes.

We do not expect the hyperbolic formulation presented herein to be preferable in every
setting, as the relaxation parameter introduces an additional tradeoff between
faithfulness to the original NLS dynamics, numerical stiffness, and flexibility of
the discretization. 
However, as has been demonstrated in the numerical experiments, 
the fact that the system is amenable to schemes tailored to hyperbolic problems
makes it attractive in a number of interesting scenarios.
Moreover, the interesting dynamics of NLSH itself, suggest that the system is worthy of further study,
both analytically and numerically.
 
Several questions remain open, and we mention those we consider most interesting.
\begin{itemize}
    \item \emph{Virial identity and blow-up.}  Whether NLSH possesses a virial-like identity
    is still unresolved.  This is closely tied to the
    question of whether the
    relaxation acts as a regularizing mechanism that can prevent finite-time blow-up
    in the focusing case.
    \item \emph{Breathers.}  The Fermi--Pasta--Ulam--Tsingou-like recurrence
    observed in Section \ref{sec:MI numerical} suggests that NLSH admits
    Akhmediev-type breather solutions, and that their amplitude and period depend
    on the relaxation parameter in a systematic way.
    \item \emph{Potentials and generalizations.} The study of problems
    with more elaborate potentials, as well as
    extensions to vector NLS and Gross--Pitaevskii systems with rotation, are
    natural next steps.
    \item \emph{Flexible and efficient discretization.} 
    From the numerical perspective, a discretization of NLSH
    that shares the flexibility of the finite volume scheme presented herein, but is more efficient
    than the fractional stepping approach, is desirable.
    To evolve only the algebraic source term implicitly in an unsplit manner
    is a clear direction compatible with the machinery used in this work.
    Moreover, an extension of the conservative discretization of Ranocha and Ketcheson \cite{ranocha2026conserving} 
    to the higher-dimensional NLSH system is also a natural next step.
\end{itemize}

\section*{Acknowledgments}
Claude Opus 5 (Anthropic) was used to refactor the Python pseudospectral solver, and to improve the plotting routines used to generate the figures in this work.
All content and results are the work of the authors, who are responsible for their accuracy and integrity.

\appendix
\section{Rotational invariance}
\label{sec:rotational-invariance}
The system \eqref{eq:ND-NLSH-real-imaginary} is invariant under rotations in 2 and 3 dimensions
in the sense of \cite[Sec. 3.2]{toro2013riemann}.
More precisely, when $N=2$, for any angle $\theta\in [0,2\pi[$,
we have
\begin{align}
    \cos(\theta) B_1 + \sin(\theta) B_2 = R_N^{-1} B_1 R_N,
\end{align}
where the rotation matrix $R_N$ is given by
\begin{align}
    R_N = \begin{bmatrix}
        1 & \mathbf{0}_{1\times 2} & 0 & \mathbf{0}_{1\times 2} \\
        \mathbf{0}_{2\times 1} & R & \mathbf{0}_{2\times 1} & \mathbf{0}_{2\times 2} \\
        0 & \mathbf{0}_{1\times 2} & 1 & \mathbf{0}_{1\times 2} \\
        \mathbf{0}_{2\times 1} & \mathbf{0}_{2\times 2} & \mathbf{0}_{2\times 1} & R
    \end{bmatrix}, \quad
    R = \begin{bmatrix}
        \cos(\theta) & \sin(\theta) \\
        -\sin(\theta) & \cos(\theta)
    \end{bmatrix}.
\end{align}
Similarly, when $N=3$, for any angles $\alpha,\beta\in [0,2\pi[$,
we have
\begin{align}
    \cos(\alpha)\cos(\beta) B_1 + \cos(\alpha)\sin(\beta) B_2 + \sin(\alpha) B_3 = R_N^{-1} B_1 R_N,
\end{align}
where 
\begin{align}
    R_N = \begin{bmatrix}
        1 & \mathbf{0}_{1\times 3} & 0 & \mathbf{0}_{1\times 3} \\
        \mathbf{0}_{3\times 1} & R & \mathbf{0}_{3\times 1} & \mathbf{0}_{3\times 3} \\
        0 & \mathbf{0}_{1\times 3} & 1 & \mathbf{0}_{1\times 3} \\
        \mathbf{0}_{3\times 1} & \mathbf{0}_{3\times 3} & \mathbf{0}_{3\times 1} & R
    \end{bmatrix}, \quad
    R = \begin{bmatrix}
        \cos(\alpha)\cos(\beta) & \cos(\alpha)\sin(\beta) & \sin(\alpha)\\
        -\sin(\beta) & \cos(\beta) & 0\\
        -\sin(\alpha)\cos(\beta) & -\sin(\alpha)\sin(\beta) & \cos(\alpha)
    \end{bmatrix}.
\end{align}
\section{A fully explicit unsplit discretization}
\label{sec:fully-explicit}
\begin{figure}
    \centering
    \includegraphics[width=0.9\textwidth]{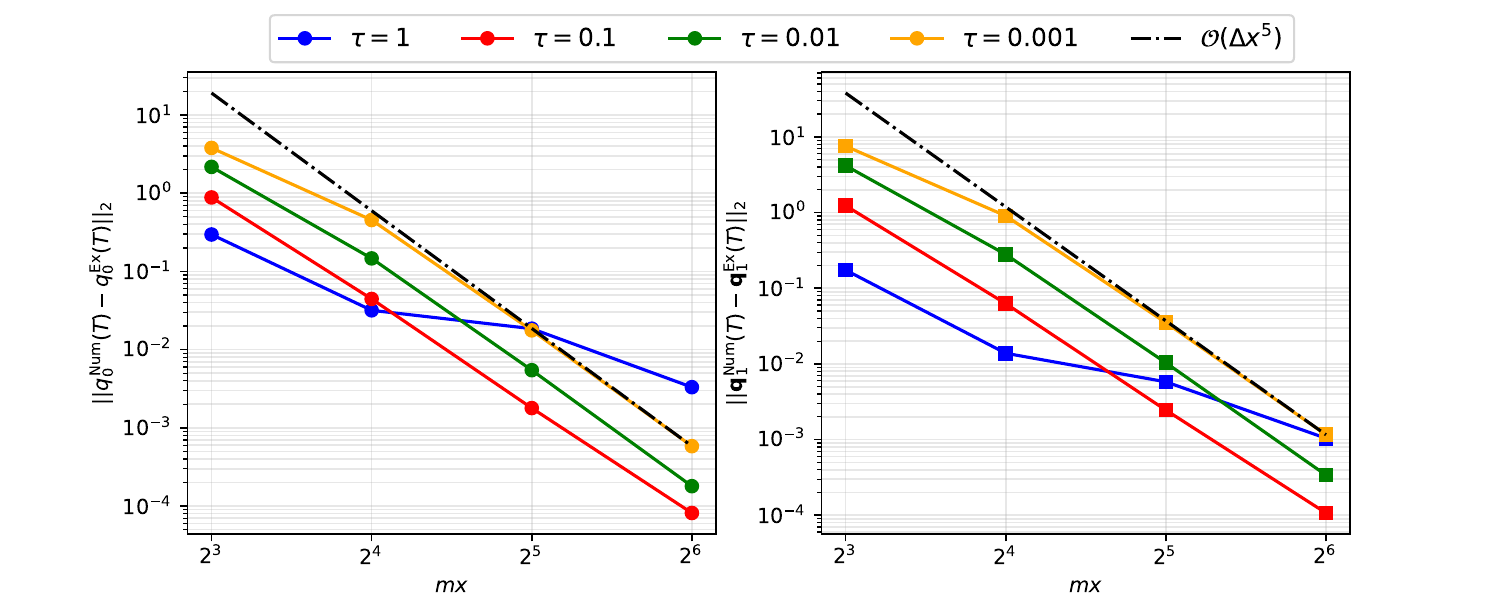}
    \caption{Convergence of the finite volume scheme for two-dimensional plane wave solutions of the NLSH system.
    The $L^2$ error is computed between the numerical solution and the exact solution for different grid resolutions
    (see Secction \ref{sec:fully-explicit} for details).}
    \label{fig: FV convergence sharpclaw}
\end{figure}

Another approach which avoids globally implicit methods for solving the NLSH system is to use a fully explicit scheme in time, evolving 
the full system \eqref{eq:ND-NLSH-real-imaginary} without splitting.
To illustrate this approach, we reproduce the test in Section \ref{sec: verification FV} using
a fifth-order WENO scheme in space and an SSPRK(3,3) scheme in time \cite{ketcheson2013high}.
Of course, for a small enough relaxation parameter $\tau$, the time step will not be limited by the CFL condition,
but rather by the stiffness of the source term, which is $\bigo(\tau^{-1})$.
Thus, we choose $\dt\propto \tau \dx$ to ensure stability of the scheme.
As can be seen in Figure \ref{fig: FV convergence sharpclaw}, 
this has also the effect of isolating the error due to the spatial discretization.
Moreover, the error constant grows approximately like $\tau^{-1/2}$ rather than 
$\tau^{-1}$.

\bibliographystyle{plain}
\bibliography{refs}

@misc{munoz2026repository,
  author       = {Muñoz-Moncayo, Carlos and Biswas, Abhijit and Ketcheson, David I. and Manvendra, Rajvanshi},
  title        = {Reproducibility repository for the manuscript "{A} hyperbolic nonlinear {S}chrödinger system in multiple spatial dimensions: structure, convergence, and numerical discretization"},
  year         = {2026},
  publisher    = {GitHub},
  journal      = {GitHub repository},
  howpublished = {\url{https://github.com/carlosmunozmoncayo/NLSH_ND}},
  note         = {Release v0}
}

@article{ranocha2026conserving,
  author   = {Ranocha, Hendrik and Ketcheson, David I.},
  title    = {Conserving Mass, Momentum, and Energy for the Benjamin–Bona–Mahony, Korteweg–de Vries, and Nonlinear Schrödinger Equations},
  journal  = {SIAM Journal on Scientific Computing},
  volume   = {48},
  number   = {4},
  pages    = {A2314-A2341},
  year     = {2026},
  doi      = {10.1137/25M183314X},
  url      = { 
              
              https://doi.org/10.1137/25M183314X
              
              
              
              },
  eprint   = { 
              
              https://doi.org/10.1137/25M183314X
              
              
              
              }
}

@article{ascher1997implicit,
  author  = {Ascher, Uri M and Ruuth, Steven J and Spiteri, Raymond J},
  journal = {Applied Numerical Mathematics},
  number  = {2},
  pages   = {151--167},
  title   = {Implicit-explicit {R}unge-{K}utta methods for time-dependent partial differential equations},
  volume  = {25},
  year    = {1997}
}

@article{larsen1983open,
  title     = {Open boundaries in short wave simulations—a new approach},
  author    = {Larsen, Jesper and Dancy, Henry},
  journal   = {Coastal engineering},
  volume    = {7},
  number    = {3},
  pages     = {285--297},
  year      = {1983},
  publisher = {Elsevier}
}

@article{petviashvili1976equation,
  author  = {Petviashvili, Vladimir I.},
  title   = {Equation of an extraordinary soliton},
  journal = {Soviet Journal of Plasma Physics},
  volume  = {2},
  pages   = {257--258},
  year    = {1976}
}

@article{ketcheson2013high,
  title     = {High-order wave propagation algorithms for hyperbolic systems},
  author    = {Ketcheson, David I and Parsani, Matteo and LeVeque, Randall J},
  journal   = {SIAM Journal on Scientific Computing},
  volume    = {35},
  number    = {1},
  pages     = {A351--A377},
  year      = {2013},
  publisher = {SIAM}
}

@article{ketcheson2019relaxation,
  title     = {Relaxation Runge--Kutta methods: Conservation and stability for inner-product norms},
  author    = {Ketcheson, David I},
  journal   = {SIAM Journal on Numerical Analysis},
  volume    = {57},
  number    = {6},
  pages     = {2850--2870},
  year      = {2019},
  publisher = {SIAM}
}

@article{berger1998amrwp,
  author   = {Berger, Marsha J. and LeVeque, Randall J.},
  title    = {Adaptive Mesh Refinement Using Wave-Propagation Algorithms for Hyperbolic Systems},
  journal  = {SIAM Journal on Numerical Analysis},
  volume   = {35},
  number   = {6},
  pages    = {2298-2316},
  year     = {1998},
  doi      = {10.1137/S0036142997315974},
  url      = { 
              
              https://doi.org/10.1137/S0036142997315974
              
              
              
              },
  eprint   = { 
              
              https://doi.org/10.1137/S0036142997315974
              
              
              
              }
}

@article{dhaouadi2019extended,
  title     = {Extended Lagrangian approach for the defocusing nonlinear Schr{\"o}dinger equation},
  author    = {Dhaouadi, Firas and Favrie, Nicolas and Gavrilyuk, Sergey},
  journal   = {Studies in Applied Mathematics},
  volume    = {142},
  number    = {3},
  pages     = {336--358},
  year      = {2019},
  publisher = {Wiley Online Library}
}

@article{antoine2013computational,
  title    = {Computational methods for the dynamics of the nonlinear Schrödinger/Gross–Pitaevskii equations},
  journal  = {Computer Physics Communications},
  volume   = {184},
  number   = {12},
  pages    = {2621-2633},
  year     = {2013},
  issn     = {0010-4655},
  doi      = {https://doi.org/10.1016/j.cpc.2013.07.012},
  url      = {https://www.sciencedirect.com/science/article/pii/S0010465513002403},
  author   = {Xavier Antoine and Weizhu Bao and Christophe Besse}
}

@article{xu2005local,
  title    = {Local discontinuous Galerkin methods for nonlinear Schrödinger equations},
  journal  = {Journal of Computational Physics},
  volume   = {205},
  number   = {1},
  pages    = {72-97},
  year     = {2005},
  issn     = {0021-9991},
  doi      = {https://doi.org/10.1016/j.jcp.2004.11.001},
  url      = {https://www.sciencedirect.com/science/article/pii/S0021999104004589},
  author   = {Yan Xu and Chi-Wang Shu}
}

@article{katsaounis2018posteriori,
  author  = {Katsaounis, Theodoros and Kyza, Irene},
  title   = {A Posteriori Error Analysis for Evolution Nonlinear Schrödinger Equations up to the Critical Exponent},
  journal = {SIAM Journal on Numerical Analysis},
  volume  = {56},
  number  = {3},
  pages   = {1405-1434},
  year    = {2018},
  doi     = {10.1137/16M1108029},
  url     = { 
             
             https://doi.org/10.1137/16M1108029
             },
  eprint  = { 
             
             https://doi.org/10.1137/16M1108029
             }
}

@article{Delfour1981FinitedifferenceSO,
  title   = {Finite-difference solutions of a non-linear Schr{\"o}dinger equation},
  author  = {Michel C. Delfour and Michel Fortin and G. Payr},
  journal = {Journal of Computational Physics},
  year    = {1981},
  volume  = {44},
  pages   = {277-288},
  url     = {https://api.semanticscholar.org/CorpusID:119748169}
}

@article{bao2003numerical,
  title     = {Numerical study of time-splitting spectral discretizations of nonlinear Schr{\"o}dinger equations in the semiclassical regimes},
  author    = {Bao, Weizhu and Jin, Shi and Markowich, Peter A},
  journal   = {SIAM Journal on Scientific Computing},
  volume    = {25},
  number    = {1},
  pages     = {27--64},
  year      = {2003},
  publisher = {SIAM}
}

@article{yuen1978relationship,
  title     = {Relationship between Benjamin--Feir instability and recurrence in the nonlinear Schr{\"o}dinger equation},
  author    = {Yuen, Henry C and Ferguson Jr, Warren E},
  journal   = {The Physics of Fluids},
  volume    = {21},
  number    = {8},
  pages     = {1275--1278},
  year      = {1978},
  publisher = {American Institute of Physics}
}

@article{akrivis2003numerical,
  author   = {Akrivis, G. D. and Dougalis, V. A. and Karakashian, O. A. and Mckinney, W. R.},
  title    = {Numerical Approximation of Blow-Up of Radially Symmetric Solutions of the Nonlinear Schrödinger Equation},
  journal  = {SIAM Journal on Scientific Computing},
  volume   = {25},
  number   = {1},
  pages    = {186-212},
  year     = {2003},
  doi      = {10.1137/S1064827597332041}
}

@book{yang2010nonlinear,
  title     = {Nonlinear waves in integrable and nonintegrable systems},
  author    = {Yang, Jianke},
  year      = {2010},
  publisher = {SIAM}
}

@article{giesselmann2025convergence,
  title   = {Convergence of hyperbolic approximations to higher-order {PDE}s for smooth solutions},
  author  = {Giesselmann, Jan and Ranocha, Hendrik},
  journal = {arXiv preprint arXiv:2508.04112},
  year    = {2025}
}

@book{cazenave2003semilinear,
  title     = {Semilinear {S}chrodinger Equations},
  author    = {Cazenave, Thierry},
  volume    = {10},
  year      = {2003},
  publisher = {American Mathematical Soc.}
}

@book{toro2013riemann,
  title     = {Riemann solvers and numerical methods for fluid dynamics: a practical introduction},
  author    = {Toro, Eleuterio F},
  year      = {2013},
  publisher = {Springer Science \& Business Media}
}

@article{ketcheson2025approximation,
  title     = {Approximation of arbitrarily high-order PDEs by first-order hyperbolic relaxation},
  author    = {Ketcheson, David I and Biswas, Abhijit},
  journal   = {Nonlinearity},
  volume    = {38},
  number    = {5},
  pages     = {055002},
  year      = {2025},
  publisher = {IOP Publishing}
}

@article{biswas2025hyperbolic,
  title   = {A Hyperbolic Approximation of the Nonlinear {S}chr\"odinger Equation},
  author  = {Biswas, Abhijit and Busaleh, Laila S and Ketcheson, David I and Mu{\~n}oz-Moncayo, Carlos and Rajvanshi, Manvendra},
  journal = {arXiv preprint arXiv:2505.21424},
  year    = {2025}
}

@article{frisch1992transition,
  author    = {Frisch, T and Pomeau, Y and Rica, S},
  title     = {Transition to dissipation in a model of superflow},
  journal   = {Physical review letters},
  volume    = {69},
  number    = {11},
  pages     = {1644},
  year      = {1992},
  publisher = {APS}
}

@article{busto2021high,
  title     = {On high order {ADER} discontinuous {G}alerkin schemes for first order hyperbolic reformulations of nonlinear dispersive systems},
  author    = {Busto, Saray and Dumbser, Michael and Escalante, Cipriano and Favrie, Nicolas and Gavrilyuk, Sergey},
  journal   = {Journal of Scientific Computing},
  volume    = {87},
  number    = {2},
  pages     = {48},
  year      = {2021},
  publisher = {Springer}
}

@article{hoefer2006dispersive,
  title     = {Dispersive and classical shock waves in Bose-Einstein condensates and gas dynamics},
  author    = {Hoefer, Mark A and Ablowitz, MJ and Coddington, I and Cornell, Eric A and Engels, P and Schweikhard, V},
  journal   = {Physical Review A—Atomic, Molecular, and Optical Physics},
  volume    = {74},
  number    = {2},
  pages     = {023623},
  year      = {2006},
  publisher = {APS}
}

@article{katsaounis2018reflection,
  title     = {On the reflection of solitons of the cubic nonlinear {S}chr{\"o}dinger equation},
  author    = {Katsaounis, Theodoros and Mitsotakis, Dimitrios},
  journal   = {Mathematical Methods in the Applied Sciences},
  volume    = {41},
  number    = {3},
  pages     = {1013--1018},
  year      = {2018},
  publisher = {Wiley Online Library}
}

@article{calhoun2008logically,
  title     = {Logically rectangular grids and finite volume methods for {PDE}s in circular and spherical domains},
  author    = {Calhoun, Donna A and Helzel, Christiane and LeVeque, Randall J},
  journal   = {SIAM review},
  volume    = {50},
  number    = {4},
  pages     = {723--752},
  year      = {2008},
  publisher = {SIAM}
}

@book{sulem2007nonlinear,
  title     = {The nonlinear {S}chr{\"o}dinger equation: self-focusing and wave collapse},
  author    = {Sulem, Catherine and Sulem, Pierre-Louis},
  volume    = {139},
  year      = {2007},
  publisher = {Springer Science \& Business Media}
}

@book{leveque2002finite,
  title     = {Finite volume methods for hyperbolic problems},
  author    = {LeVeque, Randall J},
  volume    = {31},
  year      = {2002},
  publisher = {Cambridge university press}
}

@article{mandli2016clawpack,
  title     = {Clawpack: building an open source ecosystem for solving hyperbolic PDEs},
  author    = {Mandli, Kyle T and Ahmadia, Aron J and Berger, Marsha and Calhoun, Donna and George, David L and Hadjimichael, Yiannis and Ketcheson, David I and Lemoine, Grady I and LeVeque, Randall J},
  journal   = {PeerJ Computer Science},
  volume    = {2},
  pages     = {e68},
  year      = {2016},
  publisher = {PeerJ Inc.}
}

\end{document}